\documentclass[12pt]{amsart}
\usepackage[utf8]{inputenc}
\usepackage[T1]{fontenc}
\usepackage{euscript}
\usepackage{amsfonts}
\usepackage{latexsym}
\usepackage{amscd}
\usepackage{comment}
\usepackage{tabu}
\usepackage[swedish,english]{babel}
\usepackage{lmodern}
\usepackage{amsmath, amssymb, amsthm}
\usepackage{graphicx}
\usepackage{aliascnt}
\usepackage{geometry}
\usepackage{fancyhdr}
\usepackage{euscript}
\usepackage{tikz,tkz-euclide,tikz-cd}
\usepackage{float}
\usepackage{caption}
\usepackage{subcaption}
\usepackage{csquotes}
\usepackage{enumitem}
\usepackage[colorlinks=true, linkcolor=blue, citecolor=red]{hyperref}
\usepackage{setspace}

\renewcommand{\marginpar}[1]{}
\catcode`\@=12

\def\Empty{}
\newcommand\oplabel[1]{
  \def\OpArg{#1} \ifx \OpArg\Empty {} \else
  	\label{#1}
  \fi}
		
\long\def\realfig#1#2#3#4{
\begin{figure}[htbp]
\centerline{\includegraphics[width=#4]{#2}}
\caption[#1]{#3}
\oplabel{#1}
\end{figure}}

\newcommand{\comm}[1]{}

\input{psfig}
\newtheorem{thm}{Theorem}[section]
\newtheorem*{thmH*}{Theorem A}
\newtheorem*{theoremA}{Theorem A (No Rigidity)}
\newtheorem*{theoremB}{Theorem B (Universality)}
\newtheorem{cor}[thm]{Corollary}
\newtheorem{lem}[thm]{Lemma}
\newtheorem{prop}[thm]{Proposition}

\newenvironment{pf*}[1]{\proof[#1]}{\endproof}
\usepackage{euscript}

\usepackage[OT2,OT1]{fontenc}

\newcommand{\cal}[1]{{\mathcal #1}}

\newcommand{\ignore}[1]{{}}

\theoremstyle{definition}
\newtheorem{defn}{Definition}[section]

\theoremstyle{remark}
\newtheorem{rem}{Remark}[section]

\newcommand{\diam}{\operatorname{diam}}
\newcommand{\Jac}{\operatorname{Jac}}
\newcommand{\Leb}{\rm Leb}

\newcommand{\dist}{\operatorname{dist}}

\renewcommand{\mod}{\operatorname{mod}}

\newcommand{\eps}{\epsilon}

\newcommand{\norm}[1]{\left\Vert #1 \right\Vert}
\newcommand{\absv}[1]{\left\vert #1 \right\vert}

\newcommand*{\defeq}{\mathrel{\vcenter{\baselineskip0.65ex \lineskiplimit0pt
			\hbox{$\raisebox{-0.200ex}{\scriptsize.}$}\hbox{\scriptsize.}}}
	=} 

\newcommand{\aaa}[1]{{{\mathbf{#1}}}}

\newcommand{\hol}{{\aaa H}}

\numberwithin{equation}{section}
\newcommand{\thmref}[1]{\hyperref[#1]{Theorem~\ref*{#1}}}
\newcommand{\propref}[1]{\hyperref[#1]{Proposition~\ref*{#1}}}

\newcommand{\lemref}[1]{\hyperref[#1]{Lemma~\ref*{#1}}}
\newcommand{\corref}[1]{\hyperref[#1]{Corollary~\ref*{#1}}}
\newcommand{\figref}[1]{\hyperref[#1]{Figure~\ref*{#1}}}

\newcommand{\cQ}{{\cal Q}}
\newcommand{\cJ}{{\cal J}}
\newcommand{\cA}{{\cal A}}
\newcommand{\cU}{{\cal U}}
\newcommand{\cW}{{\cal W}}

\newcommand{\cV}{{\cal V}}

\newcommand{\cI}{{\cal I}}

\newcommand{\cP}{{\cal P}}
\newcommand{\cO}{{\cal O}}

\newcommand{\cR}{{\cal R}}
\newcommand{\cRG}{{\EuScript R}}
\newcommand{\cL}{{\cal L}}

\newcommand{\cE}{{\cal E}}
\newcommand{\cS}{{\cal S}}

\newcommand{\CC}{{\Bbb C}}
\newcommand{\RR}{{\Bbb R}}
\newcommand{\TT}{{\Bbb T}}
\newcommand{\ZZ}{{\Bbb Z}}
\newcommand{\NN}{{\Bbb N}}
\newcommand{\DD}{{\Bbb D}}

\renewcommand{\AA}{{\Bbb A}}
\newcommand{\QQ}{{\Bbb Q}}
\newcommand{\bD}{{\mathbf D}}
\newcommand{\bA}{{\mathbf A}}
\newcommand{\bH}{{\mathbf H}}
\newcommand{\bK}{{\mathbf K}}

\newcommand{\bB}{{\mathbf B}}
\newcommand{\bC}{{\mathbf C}}

\newcommand{\bU}{{\mathbf U}}
\renewcommand{\bK}{{\mathbf K}}

\allowdisplaybreaks

\begin{document}

\title[Universality but no rigidity for almost commuting pairs]{Universality but no rigidity for two-dimensional perturbations of almost commuting pairs
\\[8pt]}

\author{D. Gaidashev and  D. Lilja}
\date{\today}
\maketitle

\begin{abstract}
  In this paper we consider two-dimensional dissipative maps of the annulus
  which are small perturbations of one-dimensional critical circle maps. It has
  been shown in~\cite{GaY3} that such perturbations admit an attractor which is
  a non-smooth topolgical circle - a ``critical'' circle. We study conjugacies
  of the maps that admit such attractors and show that although the maps exhibit
  universality - they approach a certain normal form when looked at small scales
  - two maps in general can not be smoothly comjugate on their critical
  attractors. This result extends the paradigm of ``universality but no
  rigidity'' in two dimensions, discovered in~\cite{dCLM}, to yet another class
  of dynamical systems.
\end{abstract}

\section{Preliminaries}
\subsection{Introduction}

Rigidity of conjugacies has always been a central issue in dynamics. In this
paper we consider conjugacies of the attractors of dissipative maps of the
annulus
$$\AA_r=\{(x,y)\in\RR^2,\;|y|<r\}/\ZZ\supset \TT$$
obtained as small perturbations of critical circle maps.

We recall that a critical circle map $f$ is a $C^3$-smooth orientation
preserving homeomorphism of the circle $\TT\equiv\RR/\ZZ$ which has a single
critical point $x_0\in\TT$ whose order $n$ is an odd integer. We will assume
that $x_0=0$, and that $n=3$. An example of a family of such maps is the {\it
  Arnold's family}
$$f_{\omega}(x)=x-{2\pi} \sin 2\pi x+\omega.$$
Each $f_{\omega}$ commutes with the unit translation, and hence it projects to a
well-defined map of the circle $\TT\equiv \RR/\ZZ$, which we denote $\hat
f_{\omega}$.

For a circle homeomorphism $f$, we will denote $\rho(f)\in\TT$ its rotation
number.
As was shown by Yoccoz in \cite{Yoc}, every critical circle map $f$ with
$\rho(f)\notin\QQ$ is topologically conjugate to the rigid rotation
\[
  R_{\rho(f)}(x)\equiv x+\rho(f)\mod \ZZ.
\]

Let $F$ be a real-analytic map 
\[
  F:\AA_r\to\AA_r.
\]
We let $$\Lambda(F)=\cap_{n\in\NN}F^n(\AA_r),$$ and refer to it as the {\it
  attractor} of $F$. In the case when $f$ is a map of the circle, we can
trivially extend it to the second coordinate, setting $F_f(x,y)=(f(x),0)$; in
this case, $\Lambda(F)=\TT$.
A natural question that can be asked in this regard is whether small
two-dimensional perturbations of a map $F_f(x,y)$ when $f$ is a critical circle
map admit attractors which are ``critical'' invariant circles: topological
circles on which the dynamics is topologically, but not smoothly, conjugate to
an irrational rotation.

The work \cite{GaY3} demonstrated that the critical circles exist in typical
families, and explained the criticality phenomenon in terms of hyperbolicity of
renormalization. Specifically, maps of the annulus with a critical circle with a
fixed irrational periodic rotation number of bounded type lie in the stable
manifold of the one-dimensional hyperbolic periodic point of renormalization.
This also makes it precise how the set of maps with a critical circle is
embedded in the typical low-parameter family: the intersection of a transverse
family with the stable manifold has a codimension equal to the dimension of the
unstable manifold of renormalization.

Of course, renormalization of critical circle maps is a classical subject, and
one of the central themes in the development of modern one-dimensional dynamics.
We refer the reader to the papers \cite{Ya3,Ya4} in which the main
renormalization conjectures, known as Lanford's Program, were proved. The
preceding historical development of the subject is described in \cite{Ya3}.

Since the space of critical circle maps is ill-suited for defining a
renormalization operator (cf~\cite{Ya3}) approaches using other spaces have been
developed from which results about critical circle maps can then be deduced.

The most common definition of renormalization of critical circle maps uses the
language of {\it commuting pairs}, as described below. This approach was
pioneered in~\cite{ORSS} and~\cite{FKS}. Analytic commuting pairs provided the
setting for proving the existence of renormalization horseshoe
attractor~\cite{dF2,dFdM2,Ya4}. However, there was a conceptual difficulty in
proving hyperbolicity in this setting, as the space of analytic commuting pairs
does not possess a natural structure of a Banach manifold.

One of the possible ways to circumvent this difficulty is the so-called
\textit{cylinder renormalization}, introduced in~\cite{Ya3}. Cylinder
renormalization acts as an analytic operator on a Banach manifold of analytic
circle maps whose domain of analyticity is an annulus in the complex cylinder.
Using this setting~\cite{Ya3} proves hyperbolicity of cylinder renormalization.
It has been also used to study analytic maps with Siegel disks~\cite{Ya5,GaY1},
and critical circle maps with non-integer critical exponents~\cite{GoY}.

Another avenue, explored in~\cite{GaY2}, is the setting \textit{almost commuting
  critical pairs}. In this approach, one drops the commutation condition for
maps in the pair, and requires the commutator to disappear only up to some
order. This is equivalent to imposing a finite number of conditions on an
appropriate Banach space of analytic circle maps, which, in turn, defines a
Banach manifold in this space. \cite{GaY2} uses this approach to give a new
proof of renormalization hyperbolicity in the setting of almost commuting pairs,
and, next, applies renormalization to small two-dimensional perturbations of
critical circle maps to prove that for every periodic rotation $\rho$ number
there is an open set $\cU_\rho$ of dissipative maps of the annulus which
contains a finite codimension manifold $\cW_\rho$ of maps with a critical circle
attractor.

In this paper we continue the study of dynamics of dissipative maps of the
annulus on their critical attractors, using the concept of the {average
  Jacobian}, and prove the following result

\begin{theoremA}
  The restriction of the dynamics of any two commuting pairs with different
  average Jacobians in $\cW_\rho$ to their critical circle attractors is not
  smoothly conjugate.
\end{theoremA}

At the same time high iterates of these maps assume a universal form.

\begin{theoremB}
  Let $Z = (A,B)\in \cW_\rho$ be a commuting pair, and $n$ be the period of
  $\rho$. Then there is an $\alpha < 1$ such that the following holds for the
  sequence of renormalizations $\cR^i Z=(A_i, B_i)$ of the pair:
  \begin{equation}
    \label{eq:2DUniversality}
    B_{k n}(x,y) =
    \begin{bmatrix}
      \xi_{kn}(x) + b^{(1 + O(\alpha^{n})) q_{k n}}f(x)y(1 + O(\alpha^{k})) \\
      x
    \end{bmatrix}
  \end{equation}
  where $\xi_{kn}(x) = \pi_1 B_{k n}(x,0))$, is a sequence of critical circle
  maps converging to a periodic orbit geometrically fast, $b$ is the average
  Jacobian of $Z$, $q_{k n}$ is the $k n$-th first return time of rotation with
  the rotation number $\rho$, and $f(x)$ is a universal function which is
  uniformly bounded away from $0$ and $\infty$ and has uniformly bounded
  derivative and distortion.
\end{theoremB}

\subsection{Commuting pairs}

We will begin by recalling the theory of one-dimensional commuting pairs and
almost-commuting pairs as well as some of the important theorem about
renormalization of almost-commuting pairs. The presentation will closely follow
the introduction from~\cite{GaY3}, previous work of the first author.

\begin{defn}
  A $C^r$-smooth (or $C^\omega$) \textit{critical commuting pair}
  $\zeta=(\eta,\xi)$ consists of two $C^r$-smooth (or $C^\omega$) orientation
  preserving interval homeomorphisms $\eta:I_\eta\to \eta(I_\eta),\;
  \xi:I_{\xi}\to \xi(I_\xi)$, where
  \begin{itemize}
  \item[(I)]{$I_\eta=[0,\xi(0)],\; I_\xi=[\eta(0),0]$;}
  \item[(II)]{Both $\eta$ and $\xi$ have homeomorphic extensions to interval
      neighborhoods of their respective domains \textit{with the same degree of
        smoothness}, that is $C^r$ (or $C^\omega$), which commute,
      $\eta\circ\xi=\xi\circ\eta$;}
  \item[(III)]{$\xi\circ\eta(0)\in I_\eta$;}
  \item[(IV)]{$\eta'(x)\ne 0\ne \xi'(y) $, for all $x\in I_\eta\setminus\{0\}$,
      and all $y\in I_\xi\setminus\{0\}$;}
  \item[(V)] each of the maps $\eta$ and $\xi$ has a cubic critical point at
    $0$:
    $$\eta'(0)=\eta''(0)=\xi'(0)=\xi''(0)=0,\text{ and }\eta'''(0)\neq 0\neq \xi'''(0) .$$
  \end{itemize}
\end{defn}

We will begin by explaining the connection between critical circle maps and
critical commuting pairs. First, let $\zeta = (\eta,\xi)$ be a critical
commuting pair and consider the interval $I = [\eta(0), \xi\circ\eta(0)]$. We
create a circle by identifying the endpoints of $I$. It can then be shown that
the map $f_{\zeta}$ given by
\[
  f_{\zeta}(x)=
  \begin{cases}
    \eta\circ\xi(x) & \text{for }x\in [\eta(0),0], \\
    \eta(x) & \text{for } x\in [0,\xi\circ\eta(0)]
  \end{cases}
\]
projects to a well-defined map of the circle up to conjugation by
diffeomorphisms of the circle. For more details see, e.g.,~\cite{GaY3}.
  
Conversely, let $f$ is a critical circle mapping with rotation number
\[
  \rho(f)=\cfrac{1}{r_0+\cfrac{1}{r_1+\cfrac{1}{r_2+\dotsb}}}.
\]
Let $p_{m}/q_{m} = [r_{0},\dots,r_{m-1}]$ and let $I_{m}$ be the arc containing
$0$ and $f^{q_{m}}(0)$ but not containing $f^{q_{m+1}}(0)$. Let $\overline{f}$
denote the lift of $f$ to the real line satisfying $\overline{f}^{\prime}(0) =
0$ and $0 < \overline{f}(0) < 1$ and similarly let $\overline{I}_{m}$ denote the
intervals in $\RR$ adjacent to $0$ that project back down to $I_{m}$. Then the
pair of maps $(\eta\vert_{\overline{I}_{m}}, \xi\vert_{\overline{I}_{m+1}})$
given by
\begin{alignat*}{2}
  & \eta \colon\thinspace \overline{I}_{m} \to \RR, \quad && x \mapsto t_{1}^{-p_{m+1}}\circ \overline{f}^{q_{m+1}}(x),\\
  & \xi \colon\thinspace \overline{I}_{m+1} \to \RR, \quad && x \mapsto t_{1}^{-p_{m}}\circ \overline{f}^{q_{m}}(x)
\end{alignat*}
where $t_{1}(x) = x+1$ denotes translation by $1$, constitutes a critical
commuting pair. For simplicity a critical commuting pair constructed in this way
from a critical circle map is denoted simply by $(f^{q_{m+1}}\vert_{I_{m}},
f^{q_{m}}\vert_{I_{m+1}})$.

\begin{defn}
  The \emph{height} $\chi(\zeta)$ of a critical commuting pair $\zeta = (\eta,
  \xi)$ is the natural number $\chi(\zeta) = r\geq 1$, if it exists, such that
  $0\in [\eta^{r+1}(\xi(0)), \eta^{r}(\xi(0))]$. If no such $r$ exists then we
  say that $\chi(\zeta) = \infty$.
\end{defn}
In fact if $\chi(\zeta) = \infty$ then $\eta\vert_{I_{\eta}}$ has a fixed point.
The height allows us to make the following definition.

\begin{defn}
  We say that a critical commuting pair $\zeta=(\eta,\xi)$ is
  \emph{renormalizable} if $\chi(\zeta)<\infty$. The \emph{renormalization} of a
  renormalizable commuting pair $\zeta=(\eta, \xi)$ is the commuting pair
  \begin{center}
    ${\cal{R}}\zeta=( \widetilde{\eta^r\circ\xi}|_{ \widetilde{I_{\xi}}},\;
    \widetilde\eta|_{\widetilde{[0,\eta^r(\xi(0))]}})$
  \end{center}
  where given $\zeta = (\eta, \xi)$ the pair $\widetilde{\zeta}$ denotes the
  pair attained by linearly rescaling $\zeta$ by the factor $\xi(0)$, i.e.,
  $\widetilde{\zeta}(x) = (\xi(0))^{-1}\zeta(\xi(0)x)$ and domains
  $\widetilde{I}_{\eta} = (\xi(0))^{-1}I_{\eta} = [0,1]$ and
  $\widetilde{I}_{\xi} = (\xi(0))^{-1}I_{\xi}$ respectively. The non-rescaled
  pair $(\eta^{r}\circ\xi\vert_{I_{\xi}}, \eta\vert_{[0, \eta^{r}(\xi(0))]})$ is
  called the \emph{pre-renormalization} and will be denoted $p\mathcal{R}\zeta$.
\end{defn}

\begin{defn}
  Let $\zeta$ be a renormalizable critical commuting pair. The \emph{rotation
    number} $\rho(\zeta) \in [0,1]$ of $\zeta$ is given by the continued
  fraction $[r_{0}, r_{1}, \dots]$ where $r_{i} = \chi(\mathcal{R}^{i}\zeta)$
  and $1/\infty$ is to be interpreted as $0$.
\end{defn}
Note that if $\chi(\mathcal{R}^{i}\zeta) = \infty$ for some $i$ then
$\rho(\zeta)$ is rational. Note also that indeed the rotation number of $\zeta$
matches the rotation number of the critical circle map $f_{\zeta}$ as
constructed above.

\begin{prop}[Proposition 1.1 from~\cite{GaY3}]
  \label{rotation number}
  The rotation number of the mapping $f_{\zeta}$ is equal to $\rho(\zeta)$.
\end{prop}

\subsection{Dynamical partitions}
\label{sec:partition}

Renormalizable critical commuting pairs induce a dynamical partition of the
interval. To define this partition, first let $\zeta = (\eta,\xi)$ be a
$k$-times renormalizable critical commuting pair and let $\bar{\zeta}_{k} =
(\bar{\eta}_{k}\vert_{I_{k}}, \bar{\xi}_{k}\vert_{J_{k}})$ denote its $k$-th
pre-renormalization. Then there are multiindices $\bar{s}_{k}$ and $\bar{t}_{k}$
such that $\bar{\eta}_{k} = \zeta^{\bar{s}_{k}}$ and $\bar{\xi}_{k} =
\zeta^{\bar{t}_{k}}$, i.e., both $\bar{\eta}_{k}$ and $\bar{\xi}_{k}$ can be
written as compositions of the original $\eta$ and $\xi$. Thus we consider the
set $\mathcal{I}$ of multiindices $\bar{s} = (a_{1}, b_{1}, \dots, a_{k},
b_{k})$ and let
\[
  \zeta^{\bar{s}} = \xi^{b_{k}}\circ
  \eta^{a_{k}}\circ\dots\circ\xi^{b_{1}}\circ\eta^{a_{1}}.
\]
We introduce an ordering $\prec$ on the set $\mathcal{I}$ as by letting
$\bar{t}\prec \bar{s}$ for $\bar{s} = (a_{1}, b_{1}, \dots, a_{k}. b_{k})$ if
$\bar{t} = (a_{1}, b_{1}, \dots, a_{m}, b_{m}, c, d)$ where $m < k$ and either
$c < a_{m+1}$ and $d = 0$ or $c = a_{m+1}$ and $d < b_{m+1}$.

\begin{defn}
  Let $\zeta$ be $k$ times renormalizable and let $\zeta_{k} =
  (\eta_{k}\vert_{I_{k}}, \xi_{k}\vert_{J_{k}})$ denote its $k$-th
  pre-renormalization, where $I_{k} = [0, \xi_{k}(0)]$ and $J_{k} = [0,
  \eta_{k}(0)]$. Let the multiindices $\bar{s}_{k}$ and $\bar{t}_{k}$ be defined
  by $\eta_{k} = \zeta^{\bar{s}_{k}}$ and $\xi_{k} = \zeta^{\bar{t}_{k}}$
  respectively. The \emph{$k$-th dynamical partition} $\mathcal{P}_{k}$ of
  $\zeta$ is the set of intervals
  \[
    \mathcal{P}_{k}= \{\zeta^{\bar{w}}(I_k)\text{ for all }\bar{w}\prec
    \bar{s}_{k}\text{ and }\zeta^{\bar{w}}(J_{k})\text{ for all }\bar{w}\prec
    \bar{t}_{k}\}.
  \]
\end{defn}

Then $\mathcal{P}_{k}$ is indeed a partition of $[\eta(0), \xi(0)]$:
\begin{itemize}
\item[(a)] $\bigcup_{H\in\mathcal{P}_{k}}H=[\eta(0),\xi(0)]$;
\item[(b)] for any two distinct elements $H_{1}$ and $H_{2}$ of $\mathcal{P}_{k}$, the interiors of $H_{1}$ and $H_{2}$ are disjoint.
\end{itemize}
\begin{figure}[H]
  \centering
  \begin{tikzpicture}
    \draw (-2,0) to (5.5,0);
    \filldraw (0,0) circle [radius = 0.5mm] node(0) {} node[below left] {$0$};
    \filldraw (-2,0) circle [radius = 0.5mm] node(eta0) {} node[below left] {$\eta(0)$};
    \filldraw (5.5,0) circle [radius = 0.5mm] node(xi0) {} node[below right] {$\xi(0)$};
    \filldraw (3.5,0) circle [radius = 0.5mm] node(etaxi0) {};
    \filldraw (1.5,0) circle [radius = 0.5mm] node(etaetaxi0) {};
    \draw[->] (0) to[out = -20, in = 200] node[below] {$\xi$} (xi0);
    \draw[->] (xi0) to[out = 160, in = 20] node[above] {$\eta$} (etaxi0);
    \draw[->] (etaxi0) to[out = 160, in = 20] node[above] {$\eta$} (etaetaxi0);

    \draw (-2,-1.5) to (1.5,-1.5);
    \filldraw (0,-1.5) circle [radius = 0.5mm] node(10) {} node[below left] {$0$};
    \filldraw (-2,-1.5) circle [radius = 0.5mm] node(1xi0) {} node[below left] {$\eta(0)$};
    \filldraw (-0.5,-1.5) circle [radius = 0.5mm] node(1etaxi0) {};
    \filldraw (1.5,-1.5) circle [radius = 0.5mm] node(1eta0) {} node[below right] {$\eta^{2}\circ\xi(0)$}; 
    \draw[->] (10) to[out = 200, in = -20] node[below] {$\eta$} (1xi0);
    \draw[->] (1xi0) to[out = 20, in = 160] node[above] {$\eta^{2}\circ\xi$} (1etaxi0);

    \draw (-2,-3) to (5.5,-3);
    \filldraw (0,-3) circle [radius = 0.5mm] node(20) {} node[below left] {$0$};
    \filldraw (-2,-3) circle [radius = 0.5mm] node(2eta0) {} node[below left] {$\eta(0)$};
    \filldraw (5.5,-3) circle [radius = 0.5mm] node(2xi0) {} node[below right] {$\xi(0)$};
    \filldraw (3.5,-3) circle [radius = 0.5mm] node(2etaxi0) {};
    \filldraw (1.5,-3) circle [radius = 0.5mm] node(2etaetaxi0) {};
    \filldraw (-0.5,-3) circle [radius = 0.5mm] node(21etaxi0) {};
    \filldraw (5,-3) circle [radius = 0.5mm];
    \filldraw (3,-3) circle [radius = 0.5mm];

    \draw[dashed] (-2,0) to (-2,-3);
    \draw[dashed] (-0.5,-1.5) to (-0.5,-3);
    \draw[dashed] (0,0) to (0,-3);
    \draw[dashed] (1.5,0) to (1.5,-3);
    \draw[dashed] (3.5,0) to (3.5,-3);
    \draw[dashed] (3.5,0) to (3.5,-3);
    \draw[dashed] (3.5,0) to (3.5,-3);
    \draw[dashed] (5.5,0) to (5.5,-3);
  \end{tikzpicture}
  \caption{The first two dynamical partitions of a pair $\zeta = (\eta, \xi)$
    with rotation number $\rho(\zeta) = [2,1,1,\dots]$. The first row shows the
    first dynamical partition of $\zeta$, the second row shows the first
    dynamical partition of the pre-renormalization $\zeta_{1}$, the third row
    shows the second dynamical partition of $\zeta$.}
  \label{fig:dynamicalPartition}
\end{figure}
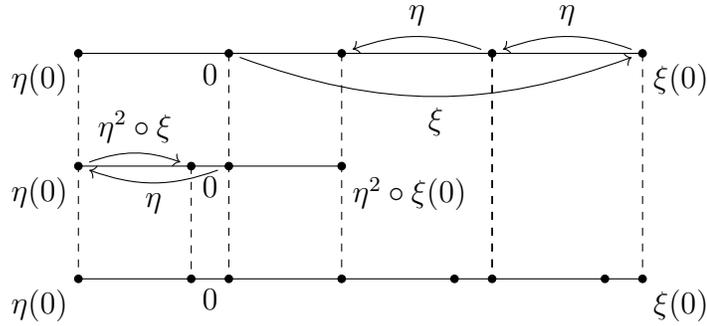

Dynamical partitions have the following property which is a consequence of the
so-called real bounds. The full statement is given in~\cite{dFdM1}. The
following reformulation is from~\cite{GaY3}.
\begin{prop}[Proposition 1.2 from~\cite{GaY3}]
  \label{realbounds}
  There exists a universal constant $C_0>1$ such that the following holds. Let
  $S$ be a compact set of $C^3$-smooth commuting pairs (note that $S$ could
  consist of a single pair). Then there exists $N=N(S)$ such that for all $n\geq
  N$ the following holds. Let $\zeta\in S$ be at least $n$ times renormalizable.
  Let $I$ and $J$ be two {\it adjacent} intervals of the $n$-th dynamical
  partition of $\zeta$. Then $I$ and $J$ are $C_0$-commensurable:
  \[
    \frac{1}{C_0}|I|<|J|<C_0|I|.
  \]
  In particular, denoting $p\cR^n\zeta=(\eta',\xi')$, we have
  \[
    \frac{1}{C_0}|I_{\xi'}|<|I_{\eta'}|<C_0|I_{\xi'}|.
  \]
\end{prop}

\subsection{Renormalization horseshoe}

The final important aspect of renormalization of critical commuting pairs we
will mention is the construction of a horseshoe attractor for renormalization of
analytic commuting pairs first presented in~\cite{Ya2}. Here $\bar{\Sigma{}}$
refers to the space of bi-infinite sequences
\[
  (\ldots,r_{-k},\ldots,r_{-1},r_0,r_1,\ldots,r_k,\ldots)\text{ with
  }r_i\in\NN\cup\{\infty\}
\]
equipped with the weak topology.

\begin{thm}[{\bf Renormalization horseshoe \cite{Ya2}}]
\label{existence of attractor}
There exists an $\mathcal{R}$-invariant set ${\cal X}$ consisting of analytic
commuting pairs with irrational rotation numbers with the following properties.
The operator $\mathcal{R}$ continuously extends to the closure (in the sense of
Carath{\'e}odory convergence)
\[
  \mathcal{A}\equiv\overline{\mathcal{X}}
\]
and the action of $\mathcal{R}$ on $\mathcal{A}$ is topologically conjugate to
the two-sided shift $\sigma:\bar{\Sigma}\to\bar{\Sigma}$:
\[
  i\circ{\mathcal{R}}\circ i^{-1}=\sigma
\]
so that if $\zeta=i^{-1}(\ldots,r_{-k},\ldots,r_{-1},r_0,r_1,\ldots,r_k,\ldots)$
then $\rho(\zeta)=[r_0,r_1,\ldots,r_k,\ldots]$. For any analytic commuting pair
$\zeta$ with an irrational rotation number we have
\[
  {\mathcal{R}}^n\zeta\to \mathcal{A}
\]
in the Carath{\'e}odory sense. Moreover, for any two analytic commuting pairs
$\zeta$, $\zeta'$ with $\rho(\zeta)=\rho(\zeta')$ we have
\[
  \dist({\mathcal{R}}^n\zeta,{\mathcal{R}}^n\zeta')\to 0
\]
for the uniform distance between analytic extensions of the renormalized pairs
on compact sets.
\end{thm}

For more information the reader is referred to~\cite{McM1} where
Carath{\'e}odory convergence is introduced in a renormalization context
and~\cite{Ya3} which contains a detailed discussion of Carath{\'e}odory
convergence for the space of analytic commuting pairs. In particular it
introduces the corresponding topology on this space.

\subsection{Spaces of analytic almost commuting pairs}

The commutation condition imposes an infinite number of conditions on the class
of commuting pairs, which makes it impossible to define a manifold structure on
this space.

 

\cite{GaY2} remedies this by introducing a class of pairs where the commutation
conditions is satisfied only approximately, and, further, proves that a space of
pairs with an appropriate commutation condition up to a finite order constitutes
a Banach manifold.

Approximate commutation has been also used by Mestel in a computer-assisted
proof of renormalization hyperbolicity by Mestel \cite{Mes}, however, \cite{Mes}
does not prove that this new class carries a manifold structure.

The following definition has been introduced in \cite{GaY2}.
\begin{defn}\label{def2}
  The space $\bB$ of {\it analytic alsmost commuting pairs} consists of pairs of
  non-decreasing interval maps
  \[
    \eta:[0,\xi(0)]\to[\eta(0),\eta\circ\xi(0)],\;\xi:[\eta(0),0]\to[\xi\circ\eta(0),\xi(0)]
  \]
  which have the following properties:
  \begin{enumerate}
  \item there exists an open neighborhood of the interval $[0,\xi(0)]$ on which
    the map $\eta$ is analytic, with a single critical point of order $3$ at the
    origin;
  \item similarly, there exists an open neighborhood of the interval
    $[\eta(0),0]$ on which the map $\xi$ is analytic, with a single critical
    point of order $3$ at the origin;
  \item the commutator
    \[
      [\eta,\xi](x)\equiv \eta\circ\xi(x)-\xi\circ\eta(x)=o(x^3)\text{ at }x=0.
    \]
  \end{enumerate}
\end{defn}

The following has been demonstrated in \cite{GaY2}.
\begin{prop}
  \label{invariance1}
  The space $\bB$ is renormalization invariant: let $\zeta\in \bB$ and
  $\rho(\zeta)\neq 0$. Then $\cR(\zeta)\in \bB$. Moreover, let
  $\rho(\zeta)\notin \QQ$. Then
  \[
    \cR^n(\zeta)\to\cA
  \]
  at a geometric rate, where $\cA$ is the hyperbolic horseshoe attractor of
  renormalization constructed in \hyperref[existence of
  attractor]{Theorem~\ref*{existence of attractor}}.
\end{prop}

Suppose, $B$ is a complex Banach space whose elements are functions of a complex
variable. Let us say that the {\it real slice} of $B$ is the real Banach space
$B^\RR$ consisting of the real-symmetric elements of $B$. If $X$ is a Banach
manifold modeled on $B$ with the atlas $\{\Psi_\gamma\}$ we shall say that $X$
is {\it real-symmetric} if
$\Psi_{\gamma_1}\circ\Psi_{\gamma_2}^{-1}(B^\RR)\subset B^\RR$ for any pair of
indices $\gamma_1$, $\gamma_2$. The {\it real slice of $X$} is then defined as
the real Banach manifold $X^\RR\subset X$ given by $\Psi_\gamma^{-1}(B^\RR)$ in
a local chart $\Psi_\gamma$. An operator $A$ defined on a subset of $X$ is {\it
  real-symmetric} if $A(X^\RR)\subset X^\RR$.

We will now restrict the definition of the space of almost commuting pairs to a
specific choice of domains of analyticity mentioned in $(1)$ and $(2)$ of the
definition \ref{def2} of $\bB$.

\begin{defn}\label{BDE}
  For a choice of topological disks $D\supset [0,1]$, $E$, we let $\bB^{D,E}_0$
  consists of pairs in $\bB$ whose maps $\eta$ and $\xi$ have bounded analytic
  continuations to $D$ and $E$ correspondingly, such that $[\eta(0),0]\subset
  E$, and such that $0$ is the only critical point of $\eta$ and $\xi$ on a
  neighborhood of $I_\eta$, $I_\xi$ respectively. We view it as a subset of the
  real slice of the complex Banach space $C^\omega(D)\times C^\omega(E)$ where
  $C^\omega(W)$ denotes the space of bounded holomorphic functions on $W$ with
  the uniform norm. Finally, denote $\bB^{D,E}$ the space of pairs in
  $\bB^{D,E}_0$ with further normalization conditions $\xi(0)=1$, and
  $\frac{1}{2C_0}<|\eta(0)|<2C_0,$ where $C_0$ is as in \propref{realbounds}.
\end{defn}

As shown in \cite{GaY2}, this space admits a manifold structure.
\begin{prop}
  \label{Banach1}
  The space $\bB^{D,E}$ equipped with the uniform norm is a real Banach
  manifold, modeled on a finite-codimensional subspace of the real slice of the
  Banach space $C^\omega(D)\times C^\omega(E)$.
\end{prop}

\subsection{Hyperbolicity and compactness of renormalization in one dimension}

We will now describe two important results that have been proved in \cite{GaY2}.

The first one is a consequence of the {\it complex a-priori bounds} and states
that renormalization improves analyticity of almost commuting pairs.

\begin{thm}\label{invariance2}
  There exists a space $\bB^{D,E}$ and $m\in\NN$ such that the following holds.
  Let $\zeta\in \bB^{D,E}$ be an $m$-times renormalizable almost commuting pair.
  There exist larger domains $D'\Supset D$, and $E'\Supset E$ so that
  \[
    \cR^m(\zeta)\in  \bB^{D',E'}.
  \]
\end{thm}

The second is the central result about the renormalization hyperbolcity.

\begin{thm}\label{th:main1}
  Let us fix a periodic point $\zeta_*\in\cA$ of $\cR$ of period $k$ and let
  $\rho_*=\rho(\zeta_*)$. There exists a space $\bB^{D,E}$ and $p=m\cdot
  k\in\NN$ such that the following holds. The pair $\zeta_*$ is a fixed point of
  $\cR^p$ in the space $\bB^{D,E}$. The image
  \[
    \cR^p(\zeta_*)\in  \bB^{D',E'} \text{ where }D'\Supset D,\;E'\Supset E.
  \]
  The linearization
  \[
    \cL\equiv D\cR^p|_{\zeta_*}
  \]
  in $\bB^{D,E}$ is a compact operator with one simple unstable eigenvalue, and
  the rest of the spectrum is compactly contained in $\DD$. The stable manifold
  $\cW^s(\zeta_*)$ of $\zeta_*$ contains all pairs in $\bB^{D,E}$ with the
  rotation number $\rho_*$.
\end{thm}

\noindent
Let $\zeta \in \cW^s(\zeta_*)$ and consider its $n$-th pre-renormalization
$\zeta_n=(\zeta^{\bar s_n},\zeta^{\bar t_n})$ defined on linear rescalings $D_n$
and $E_n$ of the sets $D$ and $E$ correspondingly. Consider the collection of
topological disks
\[
  \cV_n\equiv \{\zeta^{\bar w}(D_n)\text{ for all }\bar w\prec \bar s_n\text{ and }\zeta^{\bar w}(E_n)\text{ for all }\bar w\prec \bar t_n\}.
\]
We will refer to this collection of sets the {\it $n$-th complex dynamical
  partition} of $\zeta$. It is clear from the construction that the elements
$\zeta^{\bar w}(I_n)$ and $\zeta^{\bar w}(J_n)$ of the dynamical partition
$\cP_n$ are contained in the elements $\zeta^{\bar w}(D_n)$ and $\zeta^{\bar
  w}(E_n)$, respectively, of the complex dynamical partition $\cV_n$. Set
$\lambda_n=(-1)^n|I_n|$ so that
\[
  \cR^n\zeta(z)=\lambda_n^{-1} p\cR^n\zeta(\lambda_n z).
\]

As a consequence of \hyperref[th:main1]{Theorem \ref*{th:main1}} we have the
following:

\begin{cor}\label{cor:partition}
  Let $\zeta_*$ be as in \hyperref[th:main1]{Theorem \ref*{th:main1}}. Let
  $\zeta \in \cW^s(\zeta_*)$. Then there exists $N=N(\zeta)$, $C>0$, $C'>0$,
  $K>0$ and $0<\gamma<1$ so that for every $n>N$ the following holds.

  \begin{itemize}
  \item[$1)$] If $Q_n \in \cV_n$ then $\diam(Q_n)< C \gamma^n$.
  \item[$2)$] Any two neighboring domains $Q_n, Q'_n \in \cV_n$ are
    $K$-commensurate.
  \item[$3)$] For every $\bar w \prec \bar s_n$ (or $\bar w \prec \bar t_n$) set
    $\psi_{\bar w}^\zeta=\zeta^{\bar w} \lambda_n$. Then $\|D \psi_{\bar
      w}^\zeta |_D \|_\infty<\gamma^n$ ($\|D \psi_{\bar w}^\zeta |_E
    \|_\infty<\gamma^n$, respectively).
  \end{itemize}
\end{cor}

\section{Renormalization for dissipative two-dimensional pairs}
\label{sec:RenACM}

An extention of renormalization to a space of two-dimensional perturbations of
of one-dimensional almost commuting pairs has been been implemented in
\cite{GaY3}. \cite{GaY3} takes the route of constructing renormalization for
small perturbations of the diagonal embedding $j$ of almost commuting pairs
\[
  j((\eta,\xi))= \left(\left({\eta(x) \atop \eta(x)} \right), \left( {\xi(x)
        \atop \xi(x)} \right)  \right).
\]
Here we will follow another approach and build renormalization with a different
embedding which fixes the second component of the second map in the pair to be
equal $x$. This is an important point of departure of this paper from \cite{GaY3}. Crucially, this allows us to obtain a universality result for the second map of a high renormalization of a pair.

\subsection{Functional spaces}

Let $\zeta_*=(\eta_*,\xi_*)\in \bB^{D,E}$ be the hyperbolic fixed point of
$\cR^p$ with rotation number $\rho_*$, described in
\hyperref[th:main1]{Theorem~\ref*{th:main1}}. We denote
\[
  \bC^{D,E}\equiv (C^\omega(D)\times C^\omega(E))^\RR\supset\bB^{D,E}.
\]

We set $\Gamma_{1}=D \times \tilde{D}$, $\Gamma_{2}=E \times \tilde{E}$, where $D,
\tilde{D}, E, \tilde{E}\subset \CC$ are domains containing $0$, and let
$\bU^{\Gamma_{1}, \Gamma_{2}}$ to be the space of pairs of maps
\[
  A:\Gamma_{1}\to\CC^2,\; B:\Gamma_{2}\to \CC^2,
\]
where $A$ and $B$ are both analytic and continuous up to the boundary, equipped
with the norm
\[
  ||(A,B)||=\frac{1}{2}(||A||+||B||),\text{ where }||.||\text{ stands for the
    uniform norm.}
\]

For convenience, for a smooth function $F$ from a domain $W\subset\CC^2$ to
$\CC^2$, we will adopt the notation
\[
  ||F||_y=\underset{(x,y)\in W}{\operatorname{sup}}||\partial_yF(x,y)||.
\]

We set
\[
  \bA^{\Gamma_{1}, \Gamma_{2}}\equiv (\bU^{\Gamma_{1}, \Gamma_{2}})^\RR
\]
so that $\bA^{\Gamma_{1}, \Gamma_{2}}$ consists of pairs of real-symmetric
two-dimensional maps.

For a two-dimensional map $A(x,y)= (a(x,y),h(x,y))$ we define the projection
\[
  \cL A(x) \equiv \pi_{1}(A(x,0)) = a(x,0)
\]
Similarly, for a pair of two-dimensional maps $(A,B)(x,y)$ let us define
\[
  \cL (A,B)(x)\equiv (\pi_1(A(x,0)),\pi_1(B(x,0)))=(a(x,0),b(x,0)).
\]

Next we define the subset of all two-dimensional almost-commuting pairs, i.e.
the subset
\[
  \bD^{\Gamma_{1},\Gamma_{2}}\subset \bA^{\Gamma_{1},\Gamma_{2}}
\]
consisting of all pairs $(A,B)$ such that
\begin{align}
  \label{eq:AA} A(x,y) & = (a(x,y),h(x,y))=(a_y(x),h_y(x)),\\
  \vspace{2mm}
  \label{eq:BB} B(x,y) & = (b(x,y),x)=(b_y(x),x), 
\end{align}  
satisfying the following conditions:

\vspace{2mm}

\begin{enumerate}[label=(D\arabic*)]
\item\label{cond:D1} $\mathcal{L}[A,B](x) = \mathcal{L}(A\circ B - B\circ A)(x)
  = o\left(\absv{x}^{3}\right)$,

  \vspace{2mm}

\item\label{cond:D2} $\mathcal{L}B(0) = 1$.
\end{enumerate}

\vspace{2mm}

Let us define an embedding $\iota$ of the subset of the twice-renormalizable pairs of maps in the manifold $\bB^{D,E}$ into $\bD^{\Gamma_{1},\Gamma_{2}}$ as follows.
Consider a twice renormalizable pair $(\eta,\xi)\in \bB^{D,E}$ and let $\Lambda$ denote the rescaling, so that 
\[
  \mathcal{R}^{2}(\eta,\xi) = \Lambda((\eta^{r_0}\circ\xi)^{r_1}\circ\eta,
  \eta^{r_0}\circ\xi)
\]
where $r_0$ is the height of the pair $(\eta,\xi)$ and $r_1$ is the height
of the pair $(\eta^{r_0}\circ\xi, \eta)$. Now define
\begin{align*}
  \iota(\eta, \xi)(x,y) =
  \Lambda\left(  \begin{pmatrix}
    (\eta^{r_0}\circ\xi)^{r_1}\circ\eta(x) \\
    (\eta^{r_0}\circ\xi)^{r_1-1}\circ\eta(x)
  \end{pmatrix}, \quad
  \begin{pmatrix}
    \eta^{r_0}\circ\xi(x) \\
    x
  \end{pmatrix} \right).
\end{align*}
Here, for a pair $Z = (A,B)$ the map $\Lambda$ is defined by
\[
  \Lambda(Z) = (\Lambda^{-1}\circ A\circ \Lambda, \Lambda^{-1}\circ B\circ \Lambda)
\]
where $\Lambda (x, y) = (\lambda x, \lambda y)$ and $\lambda =\mathcal{L}(B)(0)$, i.e. $\Lambda$ is the normalization that will   ensure \hyperref[cond:D2]{Condition \ref*{cond:D2}}. For a pair $Z$ such that   $\zeta \defeq \mathcal{L}Z$ is twice renormalizable we will use $\lambda_{Z}$   to denote the rescaling constant of $\zeta$, i.e.
  \[
    \lambda_{Z} = \eta^{r_0}\circ \xi(0).
  \]

By defining the embedding in this way we have
\[
  \cL\circ \iota\equiv \cR^2.
\]

Now let $\zeta_{*}\in \bB^{D,E}$ be a periodic point of the one-dimensional
renormalization and let $Z_{*} = \iota(\zeta_{*})$. By
$\bD^{\Gamma_{1},\Gamma_{2}}_{\epsilon}(Z_{*})$ we will denote the set of pairs
$Z\in\bD^{\Gamma_{1},\Gamma_{2}}$ satisfying

\vspace{2mm}

\begin{enumerate}[label={(D\arabic*})]
  \setcounter{enumi}{2}
\item\label{cond:D3} $\mathcal{L}Z\in \bB^{D,E}$ is twice-renormalizable,

  \vspace{2mm}
  
\item\label{cond:D4} $\norm{Z-Z_{*}}\leq \epsilon$.
\end{enumerate}

\vspace{2mm}

In what follows, we will demonstrate that for each $Z_{*}$ there exists $\eps>0$
and an analytic operator $\mathcal{R}$ defined in
$\bD^{\Gamma_{1},\Gamma_{2}}_{\eps}(Z_{*})$ which has the same hyperbolic properties as
the one-dimensional version. In brief the definition consists of two steps:
\begin{enumerate}
\item Defining a pre-renormalization mimicking the one-dimensional
  renormalization in the sense that its action is close to that of $\iota\circ
  \mathcal{L}$.
\item Define a projection from two-dimensional pairs $(A,B)$ onto pairs
  $(\tilde{A}, \tilde{B})$ such that $\mathcal{L}(\tilde{A},
  \tilde{B})\in\bB^{D,E}$. This projection is not dynamical, however, crucially
  for our applications, it does not affect commuting pairs $(A,B)$.
\end{enumerate}

\subsection{Defining the pre-renormalization}

Unlike in \cite{GaY3}, where the first step was to define a sufficiently high
pre-renormalization as a pair in a neighborhood of $(\eta^{-1}(0),0)$, we will
presently consider only the second pre-renormalization of the pair.
Specifically, for a pair $Z = (A,B) \in \bD^{\Gamma_{1},\Gamma_{2}}_{\eps}(Z_{*})$ such
that $\zeta_{*}=\cL(Z_{*})$ has first two heights $r_0$ and $r_1$, we set
\begin{align*}
  p(Z) & = \left( (B\circ A^{r_0})^{r_1}\circ A, B\circ A^{r_0} \right).
\end{align*}

Next we introduce the non-linear change of coordinates $H$ defined by
\[
  H(x,y) = ((\pi_{1} A^{r_0} (\cdot, y))^{-1}(x), y).
\]
The inverse transformation $(\pi_{1} A^{r_0} (\cdot, y))^{-1}$ is well defined
for each $y$ if $\eps$ is sufficiently small, i.e. $\pi_{1} A^{r_0}(\cdot, y)$
is sufficiently close to $\eta^{r_0} = (\mathcal{L}A)^{r_0}$.

Finally we define the new pre-renormalization
\begin{align*}
  \tilde{p}(A,B) & = \left(H^{-1}\circ(B\circ A^{r_0})^{r_1}\circ A\circ H, \quad
                          H^{-1}\circ B\circ A^{r_0}\circ H\right),
\end{align*}
We remark that $\tilde{p} A$ is of the form $\tilde{p} B \circ G$.
\begin{lem}
  \label{lem:preRenormalization}
  There exists $\eps>0$ such that any $Z \in \bD^{\Gamma_{1},\Gamma_{2}}_{\eps}(Z_{*})$
  is at least twice renormalizable, and
  \[
    \norm{\Lambda(\tilde{p}Z) - \iota(\cL (Z))} = O(\varepsilon),
    \quad \norm{\tilde{p}Z}_{y} = O(\epsilon^{2}).
  \]
\end{lem}
\begin{proof}
  The Cauchy estimates for the $x$-derivatives and the fact that $\norm{Z}_{y} =
  O(\epsilon)$ on compact subsets of $\Gamma_{1}$ and $\Gamma_{2}$ imply that $a(x +
  O(\epsilon), y) = a(x, 0) + O(\epsilon)$ and $b(x + O(\epsilon), y) = b(x, 0)
  + O(\epsilon)$ on compact subsets of their domains of definition.
  
  For the second map in the pair $\tilde{p}Z$ we then get:
  \begin{align*}
    H^{-1}\circ B \circ A^{r_0} \circ H(x,y) & = H^{-1}\circ B \circ A^{r_0}( (a_{y}^{-r_0}(x), y) \\
                                             & = H^{-1}\circ B (x, \tilde{y}) \\
                                             & = H^{-1}(b_{0}(x) + O(\varepsilon), x) \\
                                             & = (a_{x}^{r_0}(b_{0}(x) + O(\varepsilon)), x).
  \end{align*}

  We would like to argue now that the composition above is well-defined om
  $\Lambda(\Gamma_{2})$.

  By definition we have that $\mathcal{L}Z = (\eta,\xi)$ is a
  twice-renormalizable almost commuting one-dimensional pair in $\bB^{D,E}$. By
  continuity, if $\eps>0$ is sufficiently small, then $\mathcal{L}Z =
  (\eta,\xi)$ is $m$-times renormalizable (since $\zeta_*$ is). By
  \hyperref[invariance2]{Theorem \ref*{invariance2}}, $\cR^m \mathcal{L}Z$ is
  defined on larger domains $D'$ and $E'$ that contain $D$ and $E$ compactly. If
  $\eps>0$ is sufficiently small, then one can choose the same domains $D'$ and
  $E'$ for all $Z$ in $\bD_\eps^{\Gamma_{1},\Gamma_{2}} (Z_*)$. Therefore, all
  the compositions that enter the definition of $\iota( \mathcal{L}Z )$ are
  defined on the larger domains $\Lambda(D')$ and $\Lambda(E')$ for all for all
  $Z$ in $\bD_\eps^{\Gamma_{1},\Gamma_{2}} (Z_*)$.

  Therefore, if $\eps>0$ is sufficiently small, the composition
  $(a_{x}^{r_0}(b_{0}(x) + O(\varepsilon))$ is well-defined on $\Lambda(D)$. We
  now have that
  \begin{align}
  \nonumber  H^{-1}\circ B \circ A^{r_0} \circ H(x,y) & = (a_{x}^{r_0}(b_{0}(x) + O(\varepsilon)), x) \\
   \nonumber  & = (a_{0}^{r_0}(b_{0}(x)) + O(\varepsilon), x) \\
   \label{approxRB}  &= (\eta^{r_0}\circ \xi(x) + O(\varepsilon), x).
  \end{align}
 
  Similarly, the compositions in the first map in the pair $\tilde{p}Z$ are well
  defined on $\Lambda(\Gamma_{1})$ and are given by:
  \begin{align}
  \nonumber   H^{-1} &\circ (B\circ A^{r_0})^{r_1}\circ A\circ H(x,y) = H^{-1}\circ (B\circ A^{r_0})^{r_1}
      \circ A(a^{-r_0}_{y}(x), y) \\
 \nonumber           & = H^{-1}\circ (B\circ A^{r_0})^{r_1}
             (a^{-(r_0-1)}_{y}(x), \tilde{y}) \\
        \nonumber    & = H^{-1}\circ(B\circ A^{r_0})^{r_1-1}
             (b_{0}(a_{0}(x)) + O(\varepsilon),
             a(x) + O(\varepsilon)) \\
         \nonumber   & = H^{-1}(b_{0}(a_0^{r_0}(\dots(b_{0}(a_{0}(x)))\dots)) + O(\varepsilon), a_{0}^{r_0}(\dots(b_{0}(a_{0}(x)))\dots) + O(\varepsilon)) \\
          \nonumber  & = (a_{0}^{r_0}(b_{0}(a_0^{r_0}(\dots(b_{0}(a_{0}(x)))\dots))) + O(\varepsilon), a_{0}^{r_0}(\dots(b_{0}(a_{0}(x)))\dots) + O(\varepsilon)) \\
        \label{approxRA}   & = ((\eta^{r_0}\circ\xi)^{r_1}\circ\eta(x) + O(\varepsilon), (\eta^{r_0}\circ\xi)^{r_1-1}\circ\eta(x) + O(\varepsilon))
  \end{align}
 
  To show that $\norm{\tilde{p}Z}_{y} = O(\epsilon^{2})$ we begin by
  noting that $\partial_{y} A(x,y) = O(\epsilon)$ and $\partial_{y} B(x,y) =
  O(\epsilon)$ on compact subsets of their domains of definition since, by
  definition, $\norm{Z - Z_{*}} \leq \epsilon$. It follows from the chain rule
  that also $\partial_{y} A^{r}(x,y) = O(\epsilon)$ for any fixed $r \geq 1$.

  Next we would like to show that $\partial_{y}H(x,y) = (O(\epsilon), 1)$. This
  follows from the fact that the first coordinate is the inverse of
  $a_{y}^{r_0}$ which has $y$-derivative of size $O(\epsilon)$ by above.

  With these two facts we get that
  \[
    \partial_y (B \circ A^{r_0} \circ H)(x,y) = \partial_y B(x,\tilde{y}) =
    \partial_2 B(x,\tilde{y}) \cdot \partial_y \tilde{y} = O(\varepsilon)\cdot
    O(\varepsilon) = O(\varepsilon^{2})
  \]
  and hence the same estimate will be true for any further compositions of maps
  with bounded $y$-derivatives, which all involved maps have by definition. This
  immediately shows that the second map in the pair $\tilde{p}Z$ has
  $y$-derivative of size $O(\epsilon^{2})$ on $\Lambda(\Gamma_{2})$. Similarly for
  the first map in the pair $\tilde{p}Z$. Thus $\norm{\tilde{p}Z}_{y} =
  O(\epsilon^{2})$ on $\Lambda(\Gamma_{1} \times \Gamma_{2})$.

\end{proof}

In light of \hyperref[lem:preRenormalization]{Lemma
  \ref*{lem:preRenormalization}} we define a preliminary renormalization
according to

\begin{equation}
  \label{eq:firstRenormalization}
  \tilde{\mathcal{R}}Z = \Lambda(\tilde{p}Z).
\end{equation}

Note that if $Z = (A,B)$ is a commutative pair then so is $(\tilde{A},
\tilde{B}) = \tilde{\mathcal{R}}Z$, hence $\tilde{\mathcal{R}}$ maps commutative
pairs in $\bD^{\Gamma_{1},\Gamma_{2}}_{\epsilon}(Z_{*})$ into $\bD^{\Gamma_{1},\Gamma_{2}}$. The
next step will be to handle non-commutative pairs in
$\bD^{\Gamma_{1},\Gamma_{2}}_{\epsilon}(Z_{*})$. For this we will create a projection
defined on
$\tilde{\mathcal{R}}\left(\bD^{\Gamma_{1},\Gamma_{2}}_{\epsilon}(Z_{*})\right)$ that
acts as the identity on commutative pairs and projects other pairs back to
$\mathbf{D}^{\Gamma_{1},\Gamma_{2}}$. Recall that the properties that need to be
satisfied are
\begin{eqnarray}
  \label{commutation} \mathcal{L}[\tilde{A},\tilde{B}](x) & = & o\left(|x|^3\right),\\
  \label{normalization} \mathcal{L}\tilde B(0) & = & 1.
\end{eqnarray}

We remark, that the pair $\tilde{\mathcal{R}}Z$ satisfies
\hyperref[normalization]{Condition (\ref*{normalization})} due to the rescaling
$\Lambda$.

Define a projection
\[
  \Pi (A,B)(x,y) = 
  \left(
    \begin{pmatrix}
      a(x,y) \\
      h(x,y)
    \end{pmatrix}, \quad
    \begin{pmatrix}
      b(x,y) + cx + dx^{2} + ex^{3} + fx^{4} \\
      x
    \end{pmatrix}
  \right)
\]
where $c,d,e,f$ are parameters.

\begin{prop}
\label{prop:2Dprojection}
There exists $\eps>0$ such that the the following holds. For every pair $Z\in
\tilde{\mathcal{R}}(\bD^{\Gamma_{1},\Gamma_{2}}_{\epsilon}(Z_{*}))$ there exists a
unique tuple $(c,d,e,f)\in \DD_{L \eps^2}(0)^{\otimes 4}$ so that
\hyperref[commutation]{Condition (\ref*{commutation})} holds for $\Pi Z$.
Furthermore, the map
\[
  \Pi\colon Z \mapsto (c,d,e,f)
\]
is analytic.
\end{prop}

\noindent
The proof is a straight-forward application of the Implicit Function Theorem but
involves several long computations. It is therefore given in
\hyperref[sec:2DprojectionProof]{Appendix \ref*{sec:2DprojectionProof}}.
\\

\noindent
We define the \emph{renormalization} of a pair $Z \in
\bD_\eps^{\Gamma_{1},\Gamma_{2}}(Z_{*})$ as
\begin{equation}
  \label{2Drenorm}
  \mathcal{R}(Z)=\Pi \tilde{\mathcal{R}}(Z).
\end{equation}
\begin{rem}\label{A=BG}
  If one denotes   $\mathcal{R}(Z)=(A_1,B_1)$, then $A_1=B_1 \circ G$ where $G$ is some analytic map from the domain of $A_1$ into that of $B_1$ whose particular form depends on the combinatorics of the maps $Z$.
\end{rem}
Suppose $\zeta_{*}$ is a one-dimensional almost commuting pair such that
$\mathcal{R}^{n}\zeta_{*} = \zeta_{*}$ where $n\geq 2$ is the smallest such
number. We can then define the two-dimensional iterated renormalization operator
through pull-backs of the composition of the operators $\mathcal{R}\colon
\mathbf{D}^{\Gamma_{1},\Gamma_{2}}_{\epsilon_{i}}(\iota(\mathcal{R}^{i}\zeta_{*}))\rightarrow
\mathbf{D}^{\Gamma_{1},\Gamma_{2}}$. This gives us the operator $\mathcal{R}^{n}\colon
\mathbf{K}_{n}(Z_{*})\rightarrow \mathbf{D}^{\Gamma_{1},\Gamma_{2}}$ where
$\mathbf{K}_{n}$ is an open subset of
$\mathbf{D}^{\Gamma_{1},\Gamma_{2}}_{\epsilon}(Z_{*})$.

With this definition the statement of
\hyperref[lem:preRenormalization]{Lemma~\ref*{lem:preRenormalization}} carries
over a fortiori to $\mathcal{R}^{n}$. We are now ready to state and prove
hyperbolicity of $\mathcal{R}^{n}$.

\begin{thm}[2D renormalization hyperbolicity]
  \label{thm:2DHyperbolicity}
  Let $\zeta_*$ be a one-dimensional renormalization fixed point of
  $\mathcal{R}^{p}$ and let $Z_{*} = \iota(\zeta_{*})$. The pair $Z_{*}$ is the
  unique fixed point of $\mathcal{R}^{n}$ for any $n=k p$, $k \in \NN$, in
  $\mathbf{K}_{n}(Z_{*})$. The linearization $D\vert_{Z_{*}}\mathcal{R}$ is a compact
  operator whose non-zero eigenvalues correspond to the tangent space of the
  embedded one-dimensional pairs, $\iota\left( \bB^{D,E} \right)$, where its
  spectrum equals that of the one-dimensional renormalization operator
  $D\vert_{\zeta_{*}}\mathcal{R}$.
\end{thm}

\begin{proof}
  First we note that
  \[
    \norm{\mathcal{R}^{n}Z}_{y} = O(\epsilon^{2^{n}})
  \]
  by \hyperref[lem:preRenormalization]{Lemma \ref*{lem:preRenormalization}} and
  hence we have that
  \[
    \lim_{n\rightarrow\infty}\norm{\mathcal{R}^{n}Z}_{y} = 0.
  \]

  Now let $Z = (A,B)$ be such that $\norm{Z}_{y} = 0$ and write
  \begin{align*}
    A(x,y) & = \left( a_{0}(x), h_{0}(x) \right), \\
    B(x,y) & = \left( b_{0}(x), x \right).
  \end{align*}
  Using $(A_{1}, B_{1}) = \mathcal{R}Z$ we then have
  \begin{align*}
    A_{1}(x, y) & = H^{-1}\circ (B\circ A^{r_0})^{r_1}\circ A\circ H(x, y) \\
                & = H^{-1}\circ (B\circ A^{r_0})^{r_1}\circ A(a_{0}^{-r_0}, y) \\
                & = H^{-1}\circ B\circ (A^{r_0}\circ B)^{r_1-1}\circ A^{r_0+1}(a_{0}^{-r_0}(x), h_{0}(x)) \\
                & = H^{-1}\circ B\circ (A^{r_0}\circ B)^{r_1-1}(a_{0}(x), h_{0}(x)) \\
                & = H^{-1}\circ B\left((a_{0}^{r_0}\circ b_{0})^{r_1-1}(a_{0}(x)), \tilde{y})\right) \\
                & = H^{-1}\left(b_{0}\circ (a_{0}^{r_0}\circ b_{0})^{r_1-1}(a_{0}(x)), (a_{0}^{r_0}\circ b_{0})^{r_1-1}(a_{0}(x))\right) \\
                & = \left( (a_{0}^{r_0}\circ b_{0})^{r_1}\circ a_{0}(x), (a_{0}^{r_0}\circ b_{0})^{r_1-1}\circ a_{0}(x) \right) \\
    B_{1}(x, y) & = H^{-1}\circ B\circ A^{r_0}\circ H(x, y) \\
                & = H^{-1}\circ B\circ A^{r_0}\left(a_{0}^{-r_0}(x), y\right) \\
                & = H^{-1}\circ B(x, h_{0}(a_{0}^{-1}(x))) \\
                & = H^{-1}(b_{0}(x), x) \\
                & = \left(a_{0}^{r_0}(b_{0}(x)), x\right)
  \end{align*}
  and we can see that $(A_{1}, B_{1}) =\iota(\mathcal{L}(Z))$ so
  $\mathcal{R}Z\in \iota\left( \bB^{D,E} \right)$. Hence we can see that
  $\mathcal{R}^{n}Z$ converges to $\iota\left( \bB^{D,E} \right)$ at rate
  $O\left( \epsilon^{2^{n}} \right)$. It follows that the spectrum of
  $\mathcal{R}$ in the complement to the tangent space of $\iota\left(
    \bB^{D,E} \right)$ at $Z_{*} = \iota(\zeta_{*})$ is trivial. Using
  hyperbolicity of renormalization of one-dimensional almost-commuting pairs as
  found in \cite{GaY3} completes the proof.
\end{proof}

\section{Attractors of dissipative maps}
\label{attractor}

As before, let $\cR^n(\zeta_*)=\zeta_*$. Fix $\rho_*\equiv
\rho(\zeta_*)\in(0,1)\setminus \QQ$. We will work in increments of $n=k p$ for
some $k \in \NN$ to be fixed later: the $n$-th iterate of the renormalization
operator will be denoted $\cRG$,

\[
  \cRG=\cR^n.
\]

Throughout this Section we will assume that $Z$ is an infinitely renormalizable
commuting pair.

By the definition of the renormalization operator
(\ref{eq:firstRenormalization}),
\begin{equation}
  \label{RGdef} \cRG Z= L_Z^{-1} \circ (Z^{\bar t_n^1}, Z^{\bar t_n^2} ) \circ L_Z,
\end{equation}
where
\begin{equation}\label{Lz}
L_Z= H_Z \circ \Lambda_Z \circ H_{\cR Z} \circ \Lambda_{\cR Z} \circ \ldots
H_{\cR^{n-1} Z} \circ \Lambda_{\cR^{n-1} Z},
\end{equation}
where the subscripts in the transformations indicate which maps are used in
these transformations. Recall that $Z$ is a commuting pair, thus the projection
$\Pi$ reduces to the identity.

Given $Z \in W^s_{\text{loc}}(Z_*)$, define the following branches of an
iterated function system on $\Gamma_1 \cup \Gamma_2$, indexed by
multi-indices $\bar w$:
\begin{equation}
  \nonumber   \psi_{\bar w}^Z \vert_{\Gamma_i} = Z^{\bar w} \circ L_Z, \quad \bar w \prec \bar t_{n}^i.
\end{equation}
The following Lemma has been proven in \cite{GaY3}.
\begin{lem}\label{lemhalf}
  There exists a neighborhood $\cS$ in $W^s_{\text{loc}}(Z_*)$ of $Z_*$ such
  that for every $\bar w \prec \bar t_n^i$ and $Z\in \cS$
  $$\|D \psi_{\bar w}^{Z} |_{\Gamma_i} \|_\infty < {1 \over 2}.$$
\end{lem}


Suppose that $Z \in \cS$. We will say that a collection of indices $\{\bar
w_m\}_{m=0}^{k-1}$ is {\it admissible on $\Gamma_{i_k}$, $i_k=1,2$,} if for all
$0 \le m \le k-1$, $\bar w_m \prec t_n^{i_{m+1}}$ whenever
\[
   \psi_{\bar w_{m}}^{\cRG^{m}  Z}(\Gamma_{i_{m+1}}) \subset \Gamma_{i_{m}},
\]
where $i_m \in \{1,2\}$.

We will denote the special collection of indices consisting only of zero
sequences by
\[
  \hat{0}_{k-1} = \{\bar{0}_{0}, \dots, \bar{0}_{k-1}\}.
\]
Note that $\psi^{Z}_{\bar{0}} = Z^{\bar{0}}\circ L_{Z} = L_{Z}$.

Given an admissible collection
$$\hat w_{k-1}=\{\bar w_0, \bar w_1, \ldots \bar w_{k-1}\},$$
consider the following {\it renormalization microscope}
\begin{equation}\label{ren_micro}
  \Psi^k_{{\hat w_{k-1}},Z}= \psi_{\bar w_0}^Z \circ \psi_{\bar w_1}^{\cRG Z} \circ \ldots \circ  \psi_{\bar w_{k-1}}^{\cRG^{(k-1)} Z}.
\end{equation}

By analogy with a dynamical partition of a one-dimensional commuting pair, the
collection
$$\cQ_{2 k n}= \left\{ \cQ_{\hat w_{k-1}}^i := \Psi^k_{{\hat w_{k-1}},Z} (\Gamma_i): \hat w_{k-1} \ {\rm is \ admissible \ on \ } \Gamma_i, i=1,2\right\}$$
will be refered to as the $2 k n$-th dynamical partition for the pair $Z$.

The following result (as well as its proof which we ommit) parallels that of
\cite{GaY3}.
\begin{lem} 
  Let $Z \in W^s_{\text{loc}}(Z_*)$. Then there exists a constant $C>0$ such  that for any $Z \in W^s_{\text{loc}}(Z_*)$ and any admissible $\hat{w}_{k-1}$
  \begin{equation}
    \label{miccontract}\|D \Psi^k_{\hat{w}_{k-1},Z}|_{\Gamma_i} \|_\infty < C {1 \over 2^k}.
  \end{equation}
\end{lem}

Set $T_a(x)\equiv x+a$, and
$$T_*\equiv (T_{\rho_*}|_{[-1,0]},T_{-1}|_{[0,\rho_*]}).$$
The following is straightforward.
\begin{lem}
  Let $\rho_*$ be a rotation number with periodic partial fraction expansion of
  period $n$. Then the commuting pair $T_*$ is a periodic point of
  renormalization
  $$\lambda_n^{-1} (T_*^{\bar u_n},T_*^{\bar v_n}) \circ \lambda_n=T_*,$$
  for some $\lambda_n$, $|\lambda_n|<1$.
\end{lem}

The intervals $[-1,0]$ and $[0,\rho_*]$ will be denoted $I$ and $J$. The
elements of the partition $\cP_{k n}$ in the orbit of the intervals $\lambda_{k
  n}(I)$ and $\lambda_{k n}(J)$ will be denoted $\cI_{k n}$ and $\cJ_{k n}$,
respectiely.

The following theorem is analogous to a result proved \cite{GaY3} using the idea
of the renormalization microscope. We will include its short proof which will
make clear our normalization for the conjugacy map $\phi$.
\begin{thm}
  \label{th:attractor}
  Let $\zeta_*=\cR^n(\zeta_*)$ be as above and let
  $$Z_*=(A_*,B_*)=\iota(\zeta_*)\in\bD^{\Gamma_1\Gamma_2}_{\epsilon_n}.$$
  Suppose $Z=(A,B)\in W^s_{\text{loc}}(Z_*)\subset \bK_n$, where $\bK_n$ is as in Theorem $\ref{thm:2DHyperbolicity}$, and suppose that maps   $A$ and $B$ commute, that is $A\circ B=B\circ A$, where defined.

  Then $Z$ has a minimal attractor $\Sigma_Z$ in $((D\cup E)\cap\RR)\times \RR$.   $\Sigma_Z$ is a Jordan arc, and the restriction $Z \vert_{\Sigma_{Z}}$ is   topologically but not smoothly conjugate to $T_{*}$.
\end{thm}

\begin{proof}
  Select a distinct point $(x_{\hat w_{k-1}}^i,y_{\hat w_{k-1}}^i)$ in each of
  the sets $\cQ^i_{\hat w_{k-1}} \in \cQ_{2 k n}$. Consider the $2 k n$-th
  dynamical partition $\cP_{2 k n}$ for the pair $T_*$ as defined in
  \hyperref[sec:partition]{Section~\ref*{sec:partition}}. Consider a
  piecewise-constant map $\phi_k$ sending the element
  $$I_{\hat w_{k-1}} = \Psi^k_{\hat w_{k-1},T_*}(I)$$
  to $(x_{\hat w_{k-1}}^i,y_{\hat w_{k-1}}^i)$, if $i=1$, or
  $$J_{\hat w_{k-1}}=\Psi^k_{\hat w_{k-1},T_*}(J)$$
  to $(x_{\hat w_{k-1}}^i,y_{\hat w_{k-1}}^i)$, if $i=2$. By
  $(\ref{miccontract})$, the diameters of the sets $\cQ^i_{\hat w_{k-1}}$
  decrease at a geometric rate. Thus, the maps $\phi_k$ converge uniformly to
  a continuous map $\phi_Z$ of the interval $[-1,\rho_*]$ which is a
  homeomorphism onto the image. Set
  $$\Sigma_Z \equiv \phi_Z([-1,\rho_*]).$$
  By construction,
  $$\phi_Z \circ T_*=Z \circ \phi_Z,$$
  and the curve $\Sigma_Z$ is the attractor for the pair $Z$. The conjugacy
  $\phi_Z$ cannot be $C^1$-smooth, since, the limiting pair $\zeta_*$ has a
  critical point at  $0$. 
\end{proof}

\section{The average Jacobian}

Let $Z \in W_{loc}^s(Z_*)$, and suppose $\phi_Z:[-1,\rho_*] \mapsto \Sigma_Z$ is
the topological conjugacy introduced in \hyperref[attractor]{Section
  \ref*{attractor}}. Since $Z|_{\Sigma_Z}$ is topologically conjugate to
$T_*|_{[-1,\rho_*]}$, it admits the unique ergodic invariant measure
$\mu_Z=\mu_{\Leb} \circ \phi_Z^{-1}$, where $\mu_{\Leb}$ is the Lebesgue measure
on $[-1,\rho_*]$. $\mu_Z$ is also the unique ergodic invariant measure for
$Z|_{\Gamma_{1} \cup \Gamma_{2}}$.

For iterated renormalizations and prerenormalizations we will use the notation
\begin{align}
  \label{Zren}  Z_{i} &= (A_{i}, B_{i}) := \cR^i Z, \\
  \label{Zreturn}  pZ_{kn} &= (pA_{kn}, pB_{kn}) = (Z^{\bar t^{1}_{kn}} ,Z^{\bar t^{2}_{kn}} ),
\end{align}
for $k \geq 1$. Here, the pre-renormalized pair $(Z^{\bar t^{1}_{kn}} ,Z^{\bar
  t^{2}_{kn}})$ is analytically conjugate to $\cRG^k Z$ by $L_Z \circ L_{\cRG Z}
\circ \ldots \circ L_{\cRG^{k-1} Z}$, similarly to \eqref{RGdef}.

We would also like to emphasize the difference between the two sets of indices
in use: on one hand, $\bar u_{k n}$ and $\bar v_{k n}$, used in
prerenormalizations for the irrational rotation $T_*$:
\begin{equation}\label{Treturn}
  \cR^{k n} T_*=\lambda_{k n}^{-1} ( T_*^{\bar u_{k n}}, T_*^{\bar v_{k n}}) \lambda_{k n},
\end{equation}
and, on the other hand, $\bar t^{1}_{k n}$ and $\bar t^{2}_{k n}$ as defined in
\eqref{Zreturn}. Indices $\bar u_{k n}$ and $\bar v_{k n}$ encode the sequence
of iterates of $T_*$ under which the {\it central intervals} $I_{ kn}=\lambda_{k
  n} (I)$ and $J_{k n} = \lambda_{k n} (J)$ return to $I_{k n} \cup J_{k n}$,
while, according to \eqref{Zreturn}, indices $\bar t^{1}_{kn}$ and $\bar
t^{2}_{kn}$ encode the returns of the {\it non-central arcs}
\begin{eqnarray}
  \nonumber I_{\hat 0_{k-1}}&=&L_z  \circ \ldots  \circ L_{\cRG^{k-1} Z} \phi_Z(I)=\Psi^k_{\hat 0_{k-1}}(\phi_Z(I)) \subset \cQ^1_{\hat 0_{k-1}},\\
  \nonumber J_{\hat 0_{k-1}} &=&L_z  \circ \ldots  \circ L_{\cRG^{k-1} Z} \phi_Z(J)=\Psi^k_{\hat 0_{k-1}}(\phi_Z(J)) \subset \cQ^2_{\hat 0_{k-1}}
\end{eqnarray}
to $I_{\hat 0_{k-1}} \cup J_{\hat 0_{k-1}}$.

Let
$$\mathfrak{w}:=|\bar w|$$
denote the $l_1$-norm of the multi-index $\bar w$.

Clearly, the return times of $T_* \vert_{I_{k n} \cup J_{k n}}$ and $Z \vert_{I_{\hat 0_{k-1}} \cup J_{\hat 0_{k-1}}}$ are the same, that is $|\bar u_{k   n}|=|\bar t^{1}_{kn}|$ and $|\bar v_{k n}|=|\bar t^{2}_{kn}|$, and are equal to the cardinalities of the partitions $\cI_{k n}$ and $\cJ_{k n}$,
respectively.

Consider the average Jacobian with respect to the measure $\mu_Z$,
\begin{equation}
  \label{avJ}
  b=\exp\left({1 \over 1 + \rho_*} {\int_{\Sigma_Z} \ln \Jac Z\vert_{\Sigma_{Z}}  \ d \mu_Z} \right).
\end{equation}

We will require a result about the dominance of the geometric convergence rates as compared to the growth of the return times.

\begin{lem}
  \label{ulambda}
  Let $\rho_*$ {have a periodic partial fraction expansion with    period $n$}. Suppose that $\bar u_{kn}$ and $\bar v_{kn}$ are such that
  $$ \lambda_{kn}^{-1} ( T_*^{\bar u_{kn}}, T_*^{\bar v_{kn}}) \lambda_{kn}=T_*.$$
  Then there exists constant $A_{m,n}=A_{m,n}(\rho_*)$, such that
\begin{eqnarray}
  \label{unorm}  |\bar u_{(k+m)n}||\lambda_{kn}^2| & < A_{m,n} \alpha^{k n}, \\
  \label{vnorm}  |\bar {v}_{(k+m)n}| |\lambda_{kn}^2| & < A_{m,n} \alpha^{k n},
\end{eqnarray}
where
$$\alpha=\sqrt{{1+\theta \over 2}}, \quad \theta=(\sqrt{5}-1)/2.$$
\end{lem}

\begin{proof}
  By a clasical result (cf \cite{McM2}), for a periodic
  $$\rho_*=[r_0,r_1,  \ldots, r_{n-1}, r_0, r_1, \ldots, r_{n-1}, \ldots],$$
  the scaling ratios are given by 
  \begin{equation}
    \label{lanbdakn}  |\lambda_{kn}|=(\theta_0 \theta_2 \cdots \theta_{n-1})^k,
  \end{equation}
  where
  $$\theta_i=[r_i, r_{i+1}, r_{i+2}, \ldots ].$$

  With the standard convention $p_{-1}=1, q_{-1}=0, p_0=0, q_0=1$, it is
  straightforward that for $m \ge 0$
  \begin{eqnarray}
    \nonumber q_{m}+p_{m} &=&   |\bar {u}_m|, \\
    \nonumber q_{m-1}+p_{m-1}     &=&   |\bar {v}_m|,
  \end{eqnarray}
  where, as before, 
  $${p_m \over q_m} =[r_0,r_1, \ldots, r_{m-1}].$$

  By the recursion formula for $q_m$,
  $${q_{m+1} \over q_{m}}=r_{m}+{q_{m-1} \over q_{m}}.$$ 
  Since at least one of the sequential numbers $q_{m-1}/q_{m-2}$ and   $q_m/q_{m-1}$ is larger than $\theta^{-1}$, the inverse golden mean (cf  \cite{Karp}, Lemma $1.22$), we have that either
  \begin{eqnarray}
    \label{case1} {q_{m+1} \over q_{m}} \le r_{m} +\theta&,&  \quad   {q_{m} \over q_{m-1}} \le r_{m-1} +1, \ {\rm or} \\
    \label{case2}  {q_{m+1} \over q_{m}} \le r_{m} +1&,& \quad   {q_{m} \over q_{m-1}} \le r_{m-1} +\theta.
  \end{eqnarray}
  We first consider the case when $kn$ is even. Then
  \begin{align}
    \nonumber  q_{(k+m)n} |\lambda_{kn}^2|&\le  \prod_{i=1}^{(k+m)n}   {q_{i} \over q_{i-1}}  \  \prod_{l=0}^{kn-1}  \theta_l^2 \\
     \nonumber   &\le  \left(\prod_{i=kn+1}^{(k+m)n} (r_{i-1}+ 1) \right)  \  \prod_{j=1}^{kn}{q_j  \over q_{j-1}} \theta_{j-1}^2 \\
    \label{prod1}  &\le   A_{m,n}  \  \prod_{j=1}^{kn} {q_j  \over q_{j-1}} \theta_{j-1}^2,    
  \end{align}
  where $A_{m,n}$ denotes the product in parenthesis and depends on $m$, $n$ and $\rho_*$.  Each pair of factors in the product is bounded as follows:
  \begin{equation}
    \nonumber     {q_j  \over q_{j-1}} \theta_{j-1}^2 {q_{j+1}  \over  q_j} \theta_j^2  \le  \max\left\{  \left(r_{j-1}+\theta  \right) \theta_{j-1}^2 \left(r_j+1  \right) \theta_j^2,  \left(r_{j-1}+1  \right) \theta_{j-1}^2 \left(r_j+\theta  \right) \theta_j^2\right\}.
  \end{equation}
  Consider the following pruduct of two consecutive factors:
  \begin{eqnarray}
    \nonumber  \left(r_j+ 1  \right) \theta_j^2 \left(r_{j+1}+ \theta  \right) \theta_{j+1}^2   & \le&  {2 \over (1+\theta_{j+1})^2 } \left( r_{j+1}+ \theta \right) \theta_{j+1}^2 \\
    \nonumber   & \le&  { 2   \left( r_{j+1}+ \theta \right)   \over \left(1+{1 \over \theta_{j+1}} \right)^2 } \\
    \nonumber   & \le&  { 2  \left( r_{j+1}+ \theta \right)   \over \left(r_{j+1}+1 \right)^2 } \\
    \label{prod2}   & \le&  { \left( 1+ \theta  \right)   \over 2 } =: \alpha^2 < 1.
  \end{eqnarray}
  Similarly for the product $\left(r_j+ \theta \right) \theta_j^2  \left(r_{j+1}+1 \right) \theta_{j+1}^2$.  Therefore,
  \begin{equation}
    \nonumber  q_{(k+m)n} |\lambda_{k n}^2|   \le   A_{m,n}  \  \alpha^{k n}.     
  \end{equation}
  If $kn$ is odd, then the extra factor in the product \eqref{prod1} can be
  included in the constant $A_{m,n}$.
  
  A similar computation holds for $p_{(k+m)n}$. The conclusion follows.
\end{proof}

We recall that $\cP_n$ denotes the dynamical partition of level $n$, while
$\cI_n$ and $\cJ_n$ denote the orbits of the central intervals $I_n$ and $J_n$
in $\cP_n$, respectively.

The second result that we will need is the following
\begin{lem}\label{proportion}
  Let $\rho_*$ have a periodic partial fraction expansion with  period $n$. Then there exist constants $C=C(\rho_*)>0$, $0<d<1$ and a universal  $\alpha<1$, independent of $\rho_*$, such that the following holds for any $l>m$, $l,m \in \NN$, and any interval $L$ of the partition $\cP_{m n}$:
  \begin{equation}
    \label{relmeasureI}  \left|{|\cI_{l n} \cap L| \over |L| } - d  \right| <  {C \alpha^{(l-m) n}}.
  \end{equation}
  Similalrly,
  \begin{equation}
    \label{relmeasureJ}  \left|{|\cJ_{l n} \cap L| \over |L| } - (1-d)  \right| <  {C \alpha^{(l-m) n}}.
  \end{equation}
\end{lem}
\begin{proof}
  The result is a consequence of the unique ergodicity of $T_*$.
  
  Consider the length of the interval $I_{k n}$,
  \begin{equation}
    |I_{k n}|=|p_{k n-1} - \rho_* q_{k n-1}| = {1 \over R_{k n} q_{k n-1} + q_{k n-2}},
  \end{equation}
  (cf \cite{Karp}, equation $(1.3)$) where
  $$R_i=r_{i-1}+[r_{i}, \ldots].$$
  In particular
  \begin{equation}
    \label{Bs}   |I_{k n}| |\bar u_{k n}|  = |I_{k n}| (q_{k n} + p_{k n})={q_{k n} + p_{k n} \over R_{k n} q_{k n-1}+q_{k n-2} }=B_{k n}.
  \end{equation}
  Let $p$ be a point in the interval $I_{k n}$. Then
  \begin{equation}
    \mathfrak{u}_{k n}^I:= \sum_{i=0}^{|\bar u_{k n}|-1} \chi_{I}(T_*^i(p)) 
  \end{equation}
  is the cardinality of $\cI_{k n} \cap I$. By the ergodicity of the irrational
  rotation
  \begin{equation}
    \label{Ds} \left|  {\mathfrak{u}^I_{k n} \over |\bar {u}_{k n}|} -  {1 \over 1+\rho_*} \right| =  \left|  {1 \over |\bar {u}_{k n}|} \sum_{i=0}^{|\bar {u}_{k n}|-1} \chi_{I}(T_*^i(p)) -  {1 \over 1+\rho_*} \int_{-1}^{\rho_*} \chi_{I} (x) d x \right| \le D_{{k n}},
  \end{equation}
  where the discrepancy $D_{k n}$ can be bounded through the  Erd\"os-Turan inequality: given an arbitrary $K$,
  \begin{equation}
    \nonumber D_{k n} \le { \log 2 \over \pi (K+1)}+{1 \over \pi |\bar {u}_{k n}|} \sum_{k=1}^K {1 \over k} \min\left(|\bar {u}_{k n}|,{1 \over \{ k \rho_*\}}\right),
  \end{equation}
where $\{ \cdot \}$ denotes the fractional part. Since $\rho_*$ is periodic, and thus, bounded, $\{k \rho_* \} > {c \over k}$, and
  \begin{equation} \label{ErdosTuran}
     D_{k n} \le {\log 2  \over \pi K}+{1 \over \pi  |\bar {u}_{k n}|} \sum_{k=1}^K {1 \over c} =  {\log 2  \over \pi K}+{K \over \pi  |\bar {u}_{k n}| c}.
  \end{equation}
  Finally, setting $K=\sqrt{|\bar {u}_{k n}|}$ and using $(\ref{Ds})$, we get that
  \begin{equation*}
     \left|  {\mathfrak{u}^I_{k n} \over |\bar {u}_{k n}|} -  {1 \over 1+\rho_*} \right|  \le {C \over \sqrt{|\bar {u}_{k n}|}}
  \end{equation*}
for some constant $C$ which depends on $c$ in $(\ref{ErdosTuran})$, and, therefore, on $\rho_*$. Similarly
  \begin{equation*}
     \left|  {\mathfrak{v}^J_{k n} \over |\bar {v}_{k n}|} -  {\rho_* \over 1+\rho_*} \right|   \le {C  \over \sqrt{|\bar {v}_{k n}}|}.
  \end{equation*}
  We, therefore, have
  \begin{equation*}
    {\left| \cI_{k n} \cap I  \right| \over |I|}={|I_{k n}| \mathfrak{u}^I_{k n} \over |I|}={|I_{k n}| |\bar {u}_{k n}|  \over 1+\rho_*} \left( 1+O \left({1 \over \sqrt{ |\bar {u}_{k n}|} } \right) \right) = {B_{k n} \over 1+\rho_*}\left( 1+O \left({1 \over \sqrt{ |\bar {u}_{k n}|} } \right)\right).
  \end{equation*}
  According to \eqref{Bs}
  \begin{align*}
    {B_{(k+1) n}  \over  B_{k n}}  &={R_{k n} q_{k n-1} + q_{k n-2}  \over R_{(k+1)n} q_{(k+1) n -1} + q_{(k+1) n-2}  }   {q_{(k+1)n}+ p_{(k+1)n} \over q_{kn} +  p_{k n}   } \\
    &={  ( r_{k n -1}+[r_{k n}, \ldots]  ) q_{k n-1} + q_{k n-2}  \over  ( r_{(k+1) n -1}+[r_{(k+1) n}, \ldots]  ) q_{(k+1) n -1} + q_{(k+1) n-2}  }   {q_{(k+1)n}+ p_{(k+1)n} \over q_{kn} +  p_{k n}   } \\
    &={ [r_{k n}, \ldots]  q_{k n-1} + q_{k n}  \over  [r_{(k+1) n}, \ldots]  q_{(k+1) n -1} + q_{(k+1) n}  }   {q_{(k+1)n}+ p_{(k+1)n} \over q_{kn} +  p_{k n}   } \\
    &={ [r_{k n}, \ldots]  q_{k n-1} + q_{k n}  \over  [r_{(k+1) n}, \ldots]  q_{(k+1) n -1} + q_{(k+1) n}  }   {1+ {p_{(k+1)n} \over q_{(k+1) n}} \over 1 +  { p_{k n} \over q_{k n}}   } \\
    &={ [r_{k n}, \ldots]  { q_{k n-1} \over  q_{k n} }+1 \over  [r_{(k+1) n}, \ldots]  { q_{(k+1) n -1}\over  q_{(k+1) n}}+1  }   {1+ {p_{(k+1)n} \over q_{(k+1) n}} \over 1 +  { p_{k n} \over q_{k n}}   } \\
    &={ [r_{k n}, \ldots]  [r_{k n-1}, r_{k n -2}, \ldots,r_0]+1 \over  [r_{(k+1) n}, \ldots]  [r_{(k+1) n-1}, r_{(k+1) n -2}, \ldots,r_0]+1  }   {1+ {p_{(k+1)n} \over q_{(k+1) n}} \over 1 +  { p_{k n} \over q_{k n}}   }.
  \end{align*}
  Notice that by the periodicity of  $\rho_*$ we have $ [r_{(k+1) n}, \ldots]= [r_{k n}, \ldots]=\rho_*$ and $[r_{k n-1}, r_{k n-2},  \ldots r_0]$ is the truncation of $[r_{(k+1) n-1}, r_{(k+1) n-2}, \ldots r_0]$  of order $k n$.
  \begin{align*}
    [r_0,r_1, \ldots, r_{(k+1)n+1}] -  [r_0,r_1, \ldots, r_{n+1}] & < C \sum_{i= k n}^{(k+1)n} {1 \over {\rm Fib}_{i} \ {\rm Fib}_{i+1}} < C  \theta^{2 k n}, \\
    \left|[r_{(k+1) n-1}, r_{(k+1) n-2}, \ldots r_0]-[r_{k n-1}, r_{k n-2}, \ldots r_0]   \right| & < {1 \over {\rm Fib}_n \ {\rm Fib}_{n+1}} < C \theta^{2 k n}, \\
  \end{align*}
where ${\rm Fib}_i$ are the Fibonacci numbers, we get
 \begin{equation}
   \nonumber   { B_{(k+1) n}  \over  B_{k n} } = {1+ \rho_*  [r_{k n-1}, r_{k n -2}, \ldots,r_0] \over 1+ \rho_*  [r_{(k+1) n-1}, r_{(k+1) n -2}, \ldots,r_0]}   {1+ {p_{(k+1)n} \over q_{(k+1) n}} \over 1 +  { p_{k n} \over q_{k n}   }}  \le (1+C \theta^{2 k n}).
  \end{equation}
Therefore, $B_{k n}$ converges as $k \rightarrow infty$, and 
$${\left| \cI_{k n} \cap I  \right| \over |I|}=d \left( 1+O \left({1 \over \sqrt{ |\bar {u}_{k n}|} } \right)\right).$$
Since $ |\bar {u}_{k n}|$ grows at least as fast as the   Fibonacci numbers, i.e. as $\theta^{-k n}$,
$${\left| \cI_{k n} \cap I  \right| \over |I|}=d \left( 1+O \left(\theta^{k n \over 2}\right)\right).$$
Claim $(\ref{relmeasureI})$ follows since $T_*$ is a renormalization fixed point of period $n$. 
The argument for the relative measure of $\cJ_{k n}$ is similar.
\end{proof}

The next result describes the Jacobians of the prerenormalizations of a pair.

\begin{lem}\label{jacob}
  Suppose $Z \in W_{loc}^s(Z_*)$. There exists $0<\alpha<1$ such that for any pair $l,m \in \NN$, $l>m$,
  \begin{eqnarray*} 
    \Jac p A_{l n} \vert_ {\cQ_{{\hat 0_{l-1}}}^1} & = & b^{\left(1+O(|\lambda_{m n}|) \right)  (1+O(\alpha^{(l-m) n})) |\bar u_{l n}|} (1 + C A_{l-m,n} \alpha^{m n})), \\
    \Jac p B_{l n} \vert_ {\cQ_{{\hat 0_{l-1}}}^2} & = & b^{\left(1+O(|\lambda_{m n}|) \right)  (1+O(\alpha^{(l-m) n})) |\bar v_{k n}|} (1 + C A_{l-m,n} \alpha^{m n})),
  \end{eqnarray*}
  and for any two points $z_1$ and $z_2$ in $\cQ_{{\hat 0_{l-1}}}^i$,
  \begin{equation*}
    \ln \left|{\Jac p Z_{l n}(z_1) \over \Jac p Z_{l n}(z_2)} \right|=O(\alpha^{l n}).
  \end{equation*}
\end{lem}
\begin{proof}
  Consider the dynamical partition $\cP_{m n}$, $ m \in \NN$, of $[-\rho_*, 1]$  for the pair $T_*$.
  
  Since $\Jac Z$ is a piecewise smooth function on $\Sigma_Z$ with one point of
  discontinuity at $0$, the standard estimates for midpoint Riemann sums
  $A_{mid}^{m n}$ over the intervals of the partition $\cP_{m n}$ give
  \begin{align*}
    \left|\ln b-A_{mid}^{m n} \right|&=\left|{1 \over 1 + \rho_*} {\int_{\Sigma_Z} \ln \Jac Z  \ d \mu_Z} -A_{mid}^{m n} \right| \\
    &\le {M_2  \over 24} \max\{|I_{m n}|^2,|J_{m n}|^2 \} \\
    &\le  C \lambda_{m n}^2,
  \end{align*}
  where $M_2$ is a bound on the absolute value of the second derivative of $\ln  \Jac A|_{\phi_Z(I)}$ and $\ln \Jac B|_{\phi_Z(J)}$.
  
  Let $m(P)$ denote the midpoint of an interval $P$. Define a piece-wise  constant funcion $F_l$ on $[-1,\rho_*]$ by settting
  $$F_l(x)=(\ln \Jac Z)\phi_Z(m(P))$$
  whenever $x$ is in an interval $P$ in $\cP_l$. $F_l \vert_P$ denotes the  restriction of $F_l$ to an interval of $P \in \cP_l$, a constant function.
 Then, for any $l>m$,
 \begin{align*}
   A_{mid}^{m n} &= {1 \over 1+\rho_*} \int_{-\rho_*}^1 F_{m n}(x) d x \\
   &= {1 \over 1 +\rho_*} \sum_{P \in \cP_{m n}} |P| \ F_{m n} \vert_P.
 \end{align*}
 Since the value of $F_{l n} |_L$ differes from that of $F_{m n}|_P$ for every  $L \subset P$ by a factor of  $(1+O(|\lambda_{m n}|)$, we get
 \begin{equation*}
   A_{mid}^{m n}= {1 \over 1 +\rho_*} \left(1+O(|\lambda_{m n}|) \right) \sum_{P \in \cP_{m n}}{|P| \over |\cI_{l n} \cap P|} \sum_{L \in \cI_{l n} \cap P}  |L| F_{l n}\vert_L.
 \end{equation*}
 By \hyperref[proportion]{Lemma $\ref*{proportion}$},
 $${|P| \over |\cI_{l n} \cap P|}={1 \over d} \left(1+O\left(\alpha^{(l-m)n}  \right)  \right),$$
 whenever $P \in \cP_{m n}$, therefore,
 \begin{equation*}
     A_{mid}^{m n} = {1 \over d (1+\rho_*)}  \left(1+O(|\lambda_{m n}|) \right) \left(1+O\left(\alpha^{(l-m)n}  \right)  \right)  \sum_{L \in \cI_{l n}} |L| F_{l n}\vert_L.
 \end{equation*}
 Since  $|L|=|I_{l n}|$  for   every \ $L \in \cI_{l n}$,
  \begin{align*}
    A_{mid}^{l n} &= {|I_{l n}| \over d  (1+\rho_*) }  \left(1+O(|\lambda_{m n}|) \right)  \left(1+O\left(\alpha^{(l-m)n}  \right)  \right)  \sum_{L \in \cI_{l n}}  (\ln \Jac Z)(\phi_Z(m(L))) \\
                  &= {|I_{l n}| |\bar u_{l n}| \over d  (1+\rho_*) } {1 \over |\bar u_{l n}|}  \left(1 \hspace{-0.5mm}  +  \hspace{-0.5mm} O(|\lambda_{m n}|) \right) \hspace{-0.7mm} \left(1  \hspace{-0.5mm}  +  \hspace{-0.5mm}  O\left(\alpha^{(l-m)n}  \right) \hspace{-0.7mm}  \right)   \hspace{-1.0mm}  \sum_{i=0}^{|\bar u_{l n}|-1}  \hspace{-1.5mm}  (\ln \Jac Z)(Z^i(\phi_Z(m(I_{l n})))) \\
                  &=  \left(1+O(|\lambda_{m n}|) \right) \left(1+O\left(\alpha^{(l-m)n}  \right)  \right)   {1 \over |\bar u_{l n}|}  \ln \Jac p A_{l n} (\phi_Z(m( I_{l n} ))),
  \end{align*} 
  and, therefore, obtain that
  \begin{equation}
    \Jac p A_{l n}(\phi_Z(m(I_{l n}))) = b^{|\bar u_{l n}|  \left(1+O(|\lambda_{m n}|) \right)  (1+O(\alpha^{(l-m) n}))} (1+O(|\bar u_{l n}| |\lambda_{m n}|^2)).
  \end{equation}
  By \hyperref[ulambda]{Lemma~\ref*{ulambda}}, $O(|u_{l n}| |\lambda_{m n}|^2)
  \rightarrow 0$ as $m n \rightarrow \infty$ geometrically under $l-m$ fixed.

  Since $\cQ^1_{\hat 0_{l-1}}$ contains a point $\phi_Z(m(I_{l n}))$, the last
  claim follows from the standard distortion bounds on $\cQ^i_{\hat 0_{l-1}}$.
\end{proof}

By the Oseledets Multiplicative Ergodic Theorem, the dynamical system
$(Z,{\Sigma_{\cRG^k Z}},\mu_Z)$ admits a set of (one or two) characteristic
exponents. In fact characteristic exponents can be computed exactly.

\begin{lem}
  Suppose $(Z,\Sigma_Z,\mu_Z)$ is a commuting pair. Then it admits two
  characteristic exponents $\chi_0=0$ and $\chi_-=\ln b$.
\end{lem}

\begin{proof}
  Let $\phi_Z$ be the topological conjugacy between $(T_*,[-1,\rho_*]$ and
  $(Z,\Sigma_Z)$. The $ k$-th renormalization $\cRG^k Z$ is analytically
  conjugate to $p \tilde \cR^{k n} Z$ on $\phi_Z(I_{k n}) \cup \phi_Z(J_{k n})$
  by $\Psi^k_{\hat 0_{k-1},Z}$. Define $\mu^{k n}_{Z}$ as $$ \mu^{k n}_{Z}:=
  (1+\rho_*) {\mu_Z \over \mu_Z(\phi_Z(I_{k n}) \cup \phi_Z(J_{k n}))}$$ on
  $\phi_Z(I_{k n}) \cup \phi_Z(J_{k n})$ and zero otherwise. Consider the
  largest Lyapunov exponent
  \begin{align*}
    \chi_0\left( p \tilde \cR^{k n} Z \vert_{\phi_Z(I_{k n}) \cup \phi_Z(J_{k n})},\mu_Z^{k n} \right) & =\chi_0(\cRG^k Z \vert_{\Sigma_{\cRG^k Z}},\mu_{\cRG^k Z}) \\
                                                                                                       & \le {1 \over 1+\rho_*} \int_{\Sigma_{\cRG^k Z}} \ln \| D \cRG^k Z \| \ d \mu_{\cRG^k Z}.
  \end{align*}
  The last integral is bounded by some constant $C$. We have
  \begin{align*}
    \min \{|\bar u_{k n}|,|\bar v_{k n}| \} \chi_0(Z \vert_{\Sigma_Z},\mu_{Z}) &\le \chi_0\left( p \tilde \cR^{k n} Z \vert_{ \phi_Z(I_{k n}) \cup \phi_Z(J_{k n})},\mu_Z^{k n}   \right) \implies \\
    \chi_0(Z \vert_{\Sigma_Z},\mu_{Z}) &\le {C \over  \min \{|\bar u_{k n}|,|\bar v_{k n}| \} },
  \end{align*}
  and $|\bar u_{k n}|$ and $|\bar v_{k n}|$ increase geometrically with $ k \in
  \NN$, therefore, $\chi_0 \le 0$.

  Now, consider $Z$ on its domain $\Gamma_{1} \cup \Gamma_{2}$, and suppose that
  $\chi_0(Z \vert_{\Sigma_{Z}},\mu_{Z})$ is negative. Then the local stable
  manifold of the derivative cocycle has dimension $2$. By the transitivity of
  $Z \vert_{\Sigma_Z}$, for any neighborhood $\cO \subset \CC^2$ of $z \in
  \Sigma_Z$ there exists a sufficiently large iterate of $Z$ which maps $\cO$
  compactly into itself. $\cO$ contains an attracting periodic orbit, which
  contradict uniqueness of the attractor in $\Gamma_{1} \cup \Gamma_{2}$. We
  conclude that $\chi_0=0$.

  The claim follows from the fact that the sum of the characteristic exponents
  is equal to $\ln b$.
\end{proof}

\section{Universality}
\label{sec:universality}

Throughout this section $Z = (A,B)$ will denote a commuting pair in the stable
manifold of a periodic point of renormalization $\mathcal{R}^{n}Z_{*} =
\EuScript{R}Z_{*} = Z_{*}$.

\begin{lem}
  The point $(1,0)$ is an attracting fixed point of the zero branch of the
  renormalization microscope of $Z_{*}$.
\end{lem}

\begin{proof}
  Denote $Z_{*} = (A_{*}, B_{*})$ and $\zeta_{*} \defeq \mathcal{L}(Z_{*}) =
  (\eta_{*},\xi_{*})$. Then we have $H(x,y) = (\eta_{*}^{-r_0}(x),y)$. Next we
  note that $\lambda_{Z_{*}} = \eta^{r_0}_{*}(\xi_{*}(0)) = \eta^{r_0}_{*}(1)$ by
  the normalization $\mathcal{L}B(0) = 1$. Thus
  \begin{align*}
    H_{Z_{*}}\circ \Lambda_{Z_{*}}(1,0) & = H_{Z_{*}}(\lambda_{Z_{*}},0) \\
                                        & = (\eta^{-r_0}_*(\lambda_{Z_{*}}),0) \\
                                        & = (\eta^{-r_0}_*(\eta^{r_0}_{*}(1)), 0) \\
                                        & = (1,0)
  \end{align*}
  so $(1,0)$ is a fixed point of the map $H_{Z_{*}}\circ \Lambda_{Z_{*}}$, hence
  it is a fixed point of $L_{Z_{*}}$. It follows that it is the fixed point of
  the zero branch of the microscope.
\end{proof}

Motivated by the above lemma and letting, as before, $T_{1}(x,y) = (x+1,y)$ be
the translation by $1$ in the $x$-coordinate we define the maps $\ell_Z$ and
$l_Z$ as follows
\begin{align*}
  {T_{1}^{-1}\circ L_{Z}\circ T_{1}(x,y)} &= (\ell_{Z,n}(x,y), \lambda_{Z,n} y), \\
  l_{Z,n}(x) &= \ell_{Z,n}(x,0),
\end{align*}
where
\begin{equation}\label{lambdaZn}
  \lambda_{Z,n} = \lambda_{Z} \cdot \lambda_{\mathcal{R} Z}\cdot\ldots\cdot\lambda_{\mathcal{R}^{n-1}Z}.
\end{equation}
Then $l_{Z,k}$ has an attracting fixed point at $x = 0$. Denote $\sigma_{Z}$ one
of the two eigenvalues of the transformation $H_{Z} \circ \Lambda_{Z}$ at
$(1,0)$ (the other being $\lambda_Z$) then
\begin{equation}
\label{prod_sigma} \prod_{i=0}^{n-1} \sigma_{\cR^i Z} = l_{Z,n}'(0).
\end{equation}
It follows by a standard argument (see, for example, Lemma $4.3$ in \cite{Yang})
that
\begin{equation}
\label{converge}\lim_{m \rightarrow \infty} \left\| {1 \over  \prod_{k=1}^{(m+1)n}  \sigma_{\cR^{k-1} Z} }  \   l_{Z,n} \circ  l_{\cRG Z,n}  \ldots l_{\cRG^m Z,n} - {1 \over  \prod_{k=0}^{n-1}  \sigma_{\cR^{k} Z_*}^{m+1} }   \ l_{Z_*,n}^{\ \circ (m+1)} \right\| = 0
\end{equation}
uniformly on compact sets, while both sequences converge to the linearizing map
$v_{Z_{*}}$ solving the equation
\begin{equation}
\label{sigma} \prod_{k=0}^{n-1} \sigma_{\cR^k Z_*}   v_{Z_{*}} = v_{Z_{*}}\circ l_{Z_*,n}.
\end{equation}
To identify the eigenvalues $\sigma_{\cR^k Z*}$, we recall $(\ref{approxRB})$,
and take the third derivative at zero of the relation $\cR Z_*^k=Z_*^{k+1}$ for
the first component of the second map in the pair. Here $Z_*^{k}$ are the pairs
in the $n$-periodic orbit $\{ Z_*^0, Z_*^1, \ldots, Z_*^{n-1}\}$ under $\cR$:
\begin{equation}
(\eta_k^{r_{k}})'(\xi_k(0)) \cdot \xi_k^{'''}(0) \cdot \lambda_{Z_*^k}^2=\xi_{k+1}^{'''}(0),
\end{equation}
where $r_k$ is the height of the pair $(\eta_k,\xi_k)$. Therefore,
\begin{equation}
\label{lprime}l_{Z_*,n}'(0)={\lambda_{Z_*,n}^3}  \prod_{i=0}^{n-1}{ \xi_{k}^{'''}(0) \over \xi_{k+1}^{'''}(0)}={\lambda_{Z_*,n}^3} { \xi_{0}^{'''}(0) \over \xi_n^{'''}(0)}={\lambda_{Z_*,n}^3}.
\end{equation}
Additionally, by \hyperref[lemhalf]{Lemma~\ref*{lemhalf}}, $\prod_{i=l+1}^{m}
\circ \ l_{\cRG^{i} Z,n}$ contracts its domains at a rate $K 2^{-m+l}$,
therefore,
\begin{equation}\label{sigma_bound}
\left\|  \prod_{l=0}^{m-1}   l_{\cRG^l Z,n}' \circ \prod_{i=l+1}^{m} \circ \  l_{\cRG^{i} Z,n} \right\| \le C  \prod_{k=0}^{n-1} |\sigma_{\cR^k Z_*}|^m  
\end{equation}
(throughout this Section $C$, $K$, $R$, $D$, $c$, $d$ etc., will denote
immaterial constants). As a consequence of hyperbolicity we have that
$T_{1}^{-1}\circ L_{\cRG^l Z}\circ T_{1}$ converges
geometrically to $( l_{Z_*,n}(x), \lambda_{Z_*,n} y)$. Furthermore, using the
notation $(\ref{Zren})$,
\begin{align}
  \nonumber   D \psi_{\bar 0}^{\cRG^l Z} & =\prod_{k=0}^{n-1} \Lambda_{Z_{l n +k}}  \left(  \left[
 \begin{array}{c c}
   \partial_x (\pi_1 A_{l n+k}^{r_{l n +k}} ( \cdot ,y))^{-1}   &  \partial_y (\pi_1 A_{l n +k }^{r_{l n +k}} (\cdot ,y))^{-1}  \\
   0 & 1
 \end{array}
 \right]  \hspace{-1mm}\circ \hspace{-0.5mm} \Lambda_{Z_{l n +k}}  \circ \right. \\
  \nonumber & \left. \hspace{77.2mm}  \circ \hspace{-2.0mm}  \prod_{i=k+1}^{n-1}\hspace{-2mm} \circ   H_{Z_{l n +i}} \hspace{-0.5mm} \circ \hspace{-0.5mm}  \Lambda_{Z_{l n +i}}  \right) \\
  \label{tilda_funcs} &  =   \prod_{k=0}^{n-1} \left[
    \begin{array}{c c}
      \partial_x (\eta_{ln+k}^{-r_{l n +k }}  \circ \lambda_{Z_{l n +k}} )  \circ l_{Z_{l n +k+1},n-k-2} +\tilde \epsilon_{Z_{l n +k}}   &  \tilde \varepsilon_{Z_{l n +k}}  \\
      0  & \lambda_{Z_{l n+k}}
    \end{array}
    \right ],
\end{align}
where $r_i$ is the height of the pair $\cL Z_i$, and $\eta_{l n +k}(x)$ is a sequence of holomorphic functions converging to the orbit of $\eta_*$ under $\cR$, and
$$\| \tilde \epsilon_{Z_{l n +k}}\| < C \eps^{2^{l n+k}} \quad {\rm and} \quad  \| \tilde \varepsilon_{Z_{l n+k}} \| <C \eps^{2^{l n +k}}.$$
Since $n$ is a fixed integer,
\begin{align}
  \nonumber D \psi_{\bar 0}^{\cRG^l Z}  & =  \left[
    \begin{array}{c c}
      \prod_{k=0}^{n-1}  l_{Z_{l n +k},0}' \circ l_{Z_{l n +k+1},n-k-2} +\epsilon_{\cRG^l Z}  & \varepsilon_{\cRG^l Z}   \\
      0  &  \lambda_{Z_{l n+k}, n}
    \end{array}
    \right ] \\
  \label{Dpsi0} & =  \left[
    \begin{array}{c c}
      l_{\cRG^l Z,n}' +\epsilon_{\cRG^l Z}  & \varepsilon_{\cRG^l Z}   \\
      0  &  \lambda_{\cRG^l Z, n}
    \end{array}
    \right ],
\end{align}
where $\lambda_{Z,n}$ has been defined in $(\ref{lambdaZn})$, and
$$\| \epsilon_{\cRG^l Z}\| < R \eps^{2^{l n}} \quad {\rm and} \quad  \| \varepsilon_{\cRG^l Z} \| <R \eps^{2^{l n}}.$$
We can now estimate the derivative of
\begin{equation}\label{Phi}
  \Phi_{Z}^m \equiv \Psi_{{\hat 0}_{m-1},Z}^m,
\end{equation}
defined in $(\ref{ren_micro})$:
\begin{eqnarray}
  \nonumber D \Phi_{\cRG^j Z}^{m} & =& \prod_{l=0}^{m-1}  \left[
    \begin{array}{c c}
      l_{\cRG^{l+j} Z,n}' + \epsilon_{\cRG^{l+j} Z}  &  \varepsilon_{\cRG^{l+j} Z} \\
      0 &  \lambda_{\cRG^{l+j}Z, n}
    \end{array}
    \right ] \circ \prod_{i=l+1}^{m-1}  \psi_{\bar 0}^{\cRG^{i+j} Z}. 
\end{eqnarray}
Assume momentarily, towards a proof by  induction, that there exist $\alpha_1<\alpha_2<1$ and $\gamma_1 < \gamma_2<1$ such that
\begin{equation}
  \label{DPsi} D \Phi_{\cRG^j Z}^{m} = \left[
    \begin{array}{c c}
      \prod_{l=0}^{m-1}   l_{\cRG^{l+j} Z, n}' \circ \prod_{i=l+1}^{m-1} \circ \  l_{\cRG^{i+j} Z,n}   + E_{m,\cRG^j Z} &  \cE_{m,\cRG^j Z}   \\
      0 &  \prod_{l=0}^{m-1}  \lambda_{\cRG^{l+j} Z, n}
    \end{array}
    \right].
\end{equation}
where
\begin{equation}
  \nonumber  c \eps^{2^{jn }}  \alpha_1^{m} \le  \| E_{m,\cRG^j Z} \|  \le C  \eps^{2^{jn }} \alpha_2^{m}, \quad   c  \eps^{2^{jn }} \gamma_1^{m} \le  \| \cE_{m,\cRG^j Z} \|  \le   C  \eps^{2^{jn }} \gamma_2^{m}.
\end{equation}

Then
\begin{equation*} \label{DPsim1}
  D \Phi_{\cRG^j Z}^{m+1} = \left[
    \begin{array}{c c}
      \prod_{l=0}^{m}  l_{\cRG^{l+j} Z, n}' \circ \prod_{i=l+1}^{m}  \circ \  l_{\cRG^{i+j} Z,n}   + E_{m+1, \cRG^j Z}   & \cE_{m+1, \cRG^j Z} \\
      0 & \prod_{l=0}^{m}  \lambda_{\cRG^{l+j} Z,n}
    \end{array}
    \right],
\end{equation*}
where the function
\begin{equation}
  \nonumber E_{m+1,\cRG^j Z} =  \eps_{\cRG^{j+m} Z} \cdot  \prod_{l=0}^{m-1}  l_{\cRG^{l+j} Z, n}' \circ \prod_{i=l+1}^{m} \circ \  l_{\cRG^{i+j} Z,n} + E_{m, \cRG^j Z} \cdot( l_{\cRG^{m+j} Z, n}'+\eps_{\cRG^{j+m} Z}  )
\end{equation}
satisfies for sufficiently large $m$ and $n$:
\begin{align}
  \nonumber \| E_{m+1,{\cRG^{j} Z}}\|  & \le R \eps^{2^{(m+j) n}} C \prod_{k=0}^{n-1} |\sigma_{\cR^k Z_*}|^{m} + C \eps^{2^{jn }} \alpha_2^{m}  \left(  D  \prod_{k=0}^{n-1} |\sigma_{\cR^k Z_*}|  +R \eps^{2^{ (m+j) n}}  \right) \\
  \label{Ele} & = C \eps^{2^{jn }} \alpha_2^{m+1} \left( R \eps^{2^{m n}}  {\prod_{k=0}^{n-1} |\sigma_{\cR^k Z_*}|^{m}  \over  \alpha_2^{m+1}} +  D { \prod_{k=0}^{n-1} |\sigma_{\cR^k Z_*}|  \over \alpha_2}  +R  {\eps^{2^{(m+j) n}} \over \alpha_2}  \right),
\end{align}
where we have used \eqref{sigma_bound}. Here, $D \approx 1$. Again, if $m$ and
$n$ are sufficiently large, then there is an
\begin{equation}\label{alpha2}
  \alpha_2 \asymp  \prod_{k=0}^{n-1} |\sigma_{\cR^k Z_*}| \asymp  |\lambda_{Z_*,n}|^3 <1
\end{equation}
for which the expression in the parenthesis is less than $1$. On the other hand
\begin{align}
  \nonumber \| E_{m+1,{\cRG^{j} Z}}\|  & \ge  c \eps^{2^{jn }} \alpha_1^{m}  d  \prod_{k=0}^{n-1} |\sigma_{\cR^k Z_*}|  -  R \eps^{2^{(m+j) n}} C \prod_{k=0}^{n-1} |\sigma_{\cR^k Z_*}|^{m} - C \eps^{2^{jn }} \alpha_2^{m} R \eps^{2^{ (m+j) n}}  \\
  \label{Ege} & =  c \eps^{2^{jn }} \alpha_1^{m+1} \left( d { \prod_{k=0}^{n-1} |\sigma_{\cR^k Z_*}|  \over \alpha_1}  - R \eps^{2^{m n}}   {C \over c}  {\prod_{k=0}^{n-1} |\sigma_{\cR^k Z_*}|^{m}  \over \alpha_1^{m+1}} -R   {C \over c}  {\eps^{2^{(m+j) n}} \over \alpha_1} { \alpha_2^{m} \over  \alpha_1^{m}}  \right).
\end{align}
Here, $d \approx 1$. Again, if $m$ and $n$ are sufficiently large, then there is an
\begin{equation}\label{alpha1}
  \alpha_1 \asymp  \prod_{k=0}^{n-1} |\sigma_{\cR^k Z_*}| \asymp  |\lambda_{Z_*,n}|^3 <1
\end{equation}
for which the expression in parenthesis is larger than $1$.
Similarly,
\begin{equation*}
  \cE_{m+1, {\cRG^{j} Z} }   = \varepsilon_{\cRG^{j+m} Z} \cdot  \left(  \prod_{l=0}^{m-1}  l_{\cRG^{l+j},n}' \circ \prod_{i=l+1}^{m} \circ \  l_{\cRG^{i+j} Z,n} + E_{m,{\cRG^{j} Z}}  \right) + \cE_{m,\cRG^{j} Z}  \cdot  \lambda_{\cRG^{m+j} Z,n},
\end{equation*}
which implies
\begin{align*}
  \| \cE_{m+1,{\cRG^{j} Z} }\| &\le K \eps^{2^{(m+j) n}} \left(  C \prod_{k=0}^{n-1} |\sigma_{\cR^k Z_*}|^{m}  +  C \alpha_2^{m+j} \right) +  C \eps^{2^{j n}} \gamma_2^m |\lambda_{\cRG^{m+j} Z,n}|   \\
  &= C \eps^{2^{j n}}  \gamma_2^{m+1} \left(  K { \eps^{2^{(m+j) n}} \over  \eps^{2^{j n}}}  \left( {\prod_{k=0}^{n-1} |\sigma_{\cR^k Z_*}|^{m}  \over \gamma_2^{m+1}} +  {\alpha_2^{j+m} \over \gamma_2^{m+1}}  \right) +  {|\lambda_{\cRG^{m+j},n}|  \over \gamma_2} \right).
\end{align*}
The expression in parenthesis is smaller than $1$ for large $m$ and $n$ and
\begin{equation} \label{gamma2}
  \gamma_2 \asymp |\lambda_{Z_{(m+j)n},n}|,  \quad  |\lambda_{Z_{(m+j)n},n}| < \gamma_2 <1
\end{equation}
We obtain that
\begin{equation} \label{cEle}
  \| \cE_{m+1,{\cRG^{j} Z} }\|   \le C \eps^{2^{j n}}  \gamma_2^{m+1}.
\end{equation}
Finally, similarly to $(\ref{Ege})$, one obtains that there exists 
\begin{equation} \label{gamma1}
  \gamma_1 \asymp |\lambda_{Z_{(m+j)n},n}|,  \quad  |\lambda_{Z_{(m+j)n},n}| > \gamma_1
\end{equation}
such that
\begin{equation} \label{cEge}
  \| \cE_{m+1,{\cRG^{j} Z} }\|   \ge c \eps^{2^{j n}}  \gamma_1^{m+1}.
\end{equation}
We are now ready to prove the main result of this section.

\begin{thm}
  \label{thm:2DUniversality}
  Let $Z = (A,B)\in W^{s}(Z _{*})$ be a commuting pair, where
  $\mathcal{R}^{n}Z_{*} = Z_{*}$. Then there is an $\alpha < 1$ such that
  \begin{equation}
    \label{eq:2DUniversality}
    B_{kn}(x,y) =
    \begin{bmatrix}
      \xi_{kn}(x) + b^{(1 + O(\alpha^{n}))\absv{\bar{v}_{kn}}}f(x)y(1 + O(\alpha^{kn})) \\
      x
    \end{bmatrix}
  \end{equation}
  where $\xi_{kn} = \mathcal{L}(B_{kn})$, $b$ is the average Jacobian of $Z$,
  $\bar{v}_{kn}$ is given by \hyperref[Treturn]{Equation \ref*{Treturn}} and
  $f$ is a universal function which is uniformly bounded away from $0$ and
  $\infty$ and has uniformly bounded derivative and distortion.
\end{thm}

\begin{proof}
  By definition we have that
  \[
    B_{kn} = (\Phi_{Z}^{k})^{-1}\circ pB_{kn}\circ \Phi_{Z}^{k}
  \]
  and therefore we have that
  \[
    \Jac B_{kn}(x,y) = \Jac pB_{kn}\left(\Phi_{Z}^{k}(x,y)\right)\frac{\Jac
      \Phi_{Z}^{k}(x,y)}{\Jac \Phi_{Z}^{k}(B_{kn}(x,y))}.
  \]
  Using \eqref{converge}, \eqref{sigma} and \eqref{DPsi} we get that
  \[
    \frac{\Jac \Phi_{Z}^{k}(x,y)}{\Jac \Phi_{Z}^{k}(B_{kn}(x,y))} \rightarrow \frac{v^{\prime}_{Z_{*}}(x)}{v^{\prime}_{Z_{*}}(\xi_{*}(x))}
  \]
  as $k\rightarrow\infty$. Denote this limit by $f(x)$.

  Writing
  \[
    B_{kn}(x,y) =
    \begin{bmatrix}
      \xi_{kn}(x) + E_{kn}(x,y) \\
      x
    \end{bmatrix}
  \]
  for some yet undetermined $E_{kn}$ we get that
  \[
    \partial_{y}E_{kn}(x,y) = \Jac B_{kn}(x,y).
  \]
  Using the above, \hyperref[jacob]{Lemma~\ref*{jacob}} with $l=k$ and a fixed
  difference $l-m$, and integrating with respect to $y$ yields the required
  form.
\end{proof}

\section{No rigidity}

We will now apply the universality result to show that any two pairs of
commuting maps on the renormalization stable manifold, whose average Jacobians
differ, can not be smoothly conjugate on their attractors. The proof is similar
to that of the pioneering work~\cite{dCLM}.

Consider the derivatives of the non-linear transformation $ H_{\cR^k \cRG^l Z}
\circ \Lambda_{\cR^k \cRG^l Z}$, $k=0, \ldots, n-1$, that enters definition
\eqref{Lz} of $L_{\cRG^l Z}$ in the case of a commuting pair $Z \in
W^s_{loc}(Z_*)$ (we will use notations $\cR^k \cRG^{l} Z=Z_{ l n+k}$
interchangeably):
\begin{align*}
  D H_{\cR^k \cRG^l Z} \circ \Lambda_{\cR^k \cRG^l Z} &= \left( \left[
    \begin{array}{c c}
      \partial_x (\pi_1 A_{l n+k}^{r_{l n +k}} ( \cdot ,y))^{-1}   &  \partial_y (\pi_1 A_{l n +k }^{r_{l n+k}} (\cdot ,y))^{-1}   \\
      0  &  1
    \end{array}
\right ] \circ  \Lambda_{Z_{l n +k}} \right)  \cdot \Lambda_{Z_{l n +k}} \\
  & =  \left[
    \begin{array}{c c}
      \partial_x (\eta_{l n+k} \circ \lambda_{Z_{l n +k}}) +\bar \delta_{Z_{l n+k}} + \bar \epsilon_{Z_{l n+k}}   &   \bar \varepsilon_{Z_{l n+k}}    \\
      0 &  \lambda_{Z_{l n+k}}
    \end{array}
    \right ]
\end{align*}
where $\sigma_{Z*}$ is as in $(\ref{sigma})$ and $\bar \epsilon_{Z_{l n
    +k}}=O(\eps^{2^{ln+k}})$, $\bar \varepsilon_{Z_{l n+k}}=O(\eps^{2^{ln+k}})$
and $\bar \delta_{Z_{l n +k}}=O(\nu^l)$, $\nu$ being the largest eigenvalue of
$D \cRG(Z_*)$ on the tangent space to $W^s_{loc}(Z_*)$.

We will momentarily concentrate on estimating the order of the function $\bar
\varepsilon_{Z_{l n+k}}$ more carefully. To that end we consider the derivative
of the transformation $\Lambda_{\cR^k \cRG^l Z}^{-1} \circ H_{\cR^k \cRG^l
  Z}^{-1}$
\begin{align}
  \nonumber   D \Lambda_{\cR^k \cRG^l Z}^{-1} \circ H_{\cR^k \cRG^l Z}^{-1}  &  =  \Lambda_{Z_{l n +k}}^{-1} \cdot \left[
    \begin{array}{c c}
      \partial_x \pi_1 A_{l n+k}^{r_{l n +k}}   &  \partial_y \pi_1 A_{l n +k }^{r_{l n +k}}  \\
      0 & 1
    \end{array}
    \right ]\\
  \label{inverseHS} & =  \Lambda_{Z_{l n +k}}^{-1} \cdot  \left[
    \begin{array}{c c}
      \partial_x \eta_*^{r_0} +\hat \delta_{Z_{l n +k}} +\hat \epsilon_{Z_{l n +k}}   &   \partial_y \pi_1 A_{l n +k }^{r_{l n +k}} \\
      0  & 1
    \end{array}
    \right ],
\end{align}
where $\hat \epsilon_{Z_{l n+k}}=O(\eps^{2^{l n +k}})$ and $\hat \delta_{Z_{l n    +k}}=O(\nu^l)$. Additionally, since maps $B$ and $A$ commute, we have
\begin{align*}
  A_{l n +k }^{r_{l n +k}}   & = \hspace{-1.0mm}  \left(\Lambda^{-1}_{l n+ k-1} \hspace{-0.6mm}  \circ \hspace{-0.6mm}   H^{-1}_{l n +k -1} \hspace{-0.6mm}   \circ \hspace{-1.0mm}  \left( B_{l n +k-1} \hspace{-0.6mm}  \circ \hspace{-0.6mm}  A_{l n +k -1}^{r_{l n +k }}  \right)^{r_{l n +k+1}} \hspace{-3.0mm}  \circ \hspace{-0.6mm}  A_{l n +k -1} \hspace{-0.6mm}  \circ \hspace{-0.6mm}  H_{l n +k -1}  \hspace{-0.6mm} \circ \hspace{-0.6mm}  \Lambda_{l n +k -1} \right)^{r_{l n +k}}  \\
  & =  \hspace{-1mm} \left(   \Lambda^{-1}_{l n+ k-1} \hspace{-0.6mm}  \circ  \hspace{-0.6mm}  H^{-1}_{l n +k -1}  \hspace{-0.6mm}  \circ \hspace{-0.6mm}  \left( B_{l n +k-1} \hspace{-0.6mm}  \circ \hspace{-0.6mm}  A_{l n +k -1}^{r_{l n +k}}  \right)^{r_{l n +k +1}-1} \hspace{-0.6mm}  \circ \hspace{-0.6mm}  A_{l n +k -1} \hspace{-0.6mm}  \circ \hspace{-0.6mm}  H_{l n +k -1} \hspace{-0.6mm}  \circ \hspace{-0.6mm}  \Lambda_{l n +k -1} \hspace{-0.4mm}  \circ \hspace{-0.6mm}  \right.\\
  &  \phantom{=}   \hspace{16.7mm} \left. \circ \ \Lambda^{-1}_{l n+ k-1} \circ H^{-1}_{l n +k -1} \circ  B_{l n +k-1} \circ A_{l n +k -1}^{r_{l n +k}}  \circ  H_{l n +k -1} \circ \Lambda_{l n +k -1}  \right)^{r_{l n +k}}\\
  & = \hspace{-0.6mm}  G_{l n +k} \circ B_{l n +k},
\end{align*}
where $G_{l n +k}$ is some analytic map defined on the range of $B_{l n +k }$
which depends on the combinatorics of $Z$. We have therefore,
$$ \partial_y \pi_1 A_{l n +k }^{r_{l n +k}} (x,y)= \nabla \pi_1 G_{l n +k} \circ B_{l n +k } (x,y) \cdot \left(  \partial_y \pi_1 B_{l n +k }(x,y)  ,   \partial_y \pi_2 B_{l n +k }(x,y) \right).$$ 
We can now use Theorem $\ref{thm:2DUniversality}$ to get
$$ \partial_y \pi_1 A_{l n +k }^{r_{l n +k}} (x,y)= \partial_x \pi_1 G_{l n +k} \circ B_{l n +k } (x,y)  b^{(1 + O(\alpha^{n}))\absv{\bar{v}_{ln}}}f(x) (1 + O(\alpha^{l})),$$ 
and since the derivatives of the analytic map $G_{l n +k}$ are uniformly bounded
on the image of the domain of $A_{l n+k}$ under $B_{l n +k}$, we get that
\begin{equation}\label{dyAr}
  \partial_y \pi_1 A_{l n +k }^{r_{l n +k}}  = O(b^{(1 + O(\alpha^{n}))\absv{\bar{v}_{ln}}}),
\end{equation}
where $b$ is the average Jacobian of the pair $Z$. Inverting the matrix
\eqref{inverseHS}, we obtain a similar result for $\tilde \varepsilon_{Z_{l
    n+k}}$ introduced in \eqref{tilda_funcs}:
\begin{equation}
  \tilde \varepsilon_{Z_{l n+k}} = O\left(b^{(1 + O(\alpha^{n}))\absv{\bar{v}_{ln}}} \right).
\end{equation}
Repeating the calculations of Section \ref{sec:universality} leading to
\eqref{DPsi}, \eqref{Ege}, \eqref{Ele}, \eqref{cEge} and \eqref{cEle}, we get
the following bounds for the functions entering the expression \eqref{DPsi} for
the derivative of $\Phi_{\cRG^j Z}^{m}$
\begin{alignat}{2}
  \label{newE} c \eps^{2^{j n}} \alpha_1^{m} \le  & \|  E_{m,\cRG^j Z}  \| && \le   C \eps^{2^{j n}} \alpha_2^{m},\\
  \label{newcE} c  b^{(1 + O(\alpha^{n}))\absv{\bar{v}_{jn}}} \gamma_1^m  \le & \| \cE_{m,\cRG^j Z}  \| && \le C  b^{(1 + O(\alpha^{n}))\absv{\bar{v}_{jn}}} \gamma_2^m 
\end{alignat}
with $\alpha_1$, $\alpha_2$, $\gamma_1$ and $\gamma_2$ as in $(\ref{alpha1})$,
$(\ref{alpha2})$, $(\ref{gamma1})$ and $(\ref{gamma2})$.

\begin{thm}
  Let $Z$ and $\tilde Z$ be two diffeomorphisms in $W^s_{loc}(Z_*)$ with average
  Jacobians $b>\tilde b$. Let $\phi$ be the homeomorphism that conjugates $Z
  \vert_{\Sigma_Z}$ to $Z \vert_{\Sigma_{\tilde{Z}}}$ normalized so that
  $\phi(\tau)=\tau$.  If  $\phi$ is $\kappa$-H\"older, then
  $$\kappa \le {1 \over 3} + {2 \over 3} { \ln b \over \ln \tilde b}.$$
\end{thm}

\begin{proof}
  Throughout the proof we will use $C$, $K$, $P$, $R$ for constants whose specific values are irrelevant to the proof.
 
  We can assume that $n$ is large, possibly after taking its integer multiple.
  For such fixed $n$, choose a large $k$ and $l=l(k)$ so that
  \begin{equation}\label{relations}
    \prod_{i=0}^{n-1}|\sigma_{\cR^i Z_*}|^{l+1} \ll \tilde  b^{(1+O(\alpha^n))|\bar \upsilon_{k n}| }  |\lambda_{Z_*,n}|^l \ll  \prod_{i=0}^{n-1}|\sigma_{\cR^i Z_*}|^{l}  \ll b^{(1+O(\alpha^n))|\bar \upsilon_{k n}| } \gamma_1^l.
  \end{equation}
  Such a choice is possible, since by \eqref{prod_sigma}, \eqref{lprime} and
  \eqref{gamma1} the first, the third and the fourth expressions in
  \eqref{relations} are commensurate with $\lambda_{Z_*,n}^{3 (l+1)}$,
  $\lambda_{Z_*,n}^{3 l}$ and $ b^{(1+O(\alpha^n))|\bar \upsilon_{k n}| }
  \lambda_{Z_*,n}^l$ respectively, and, therefore, \eqref{relations} is implied
  by
  \begin{equation}\label{relations2}
    |\lambda_{Z_*,n}|^{3 (l+1)} \ll \tilde  b^{|\bar \upsilon_{k n}| } |\lambda_{Z_*,n}|^l \ll   |\lambda_{Z_*,n}|^{3 l}  \ll b^{|\bar \upsilon_{k n}| } |\lambda_{Z_*,n}|^{l},
  \end{equation}
  which, in turn, is implied by 
  \begin{equation}\label{relations3}
    2 l+3 > |\bar \upsilon_{k n}| {\ln \tilde b \over \ln |\lambda_{Z_*,n}| } > 2 l >  |\bar \upsilon_{k n}| {\ln  b \over \ln |\lambda_{Z_*,n}| }, \quad |\lambda_{Z_*n}|  \ll 1.
  \end{equation}
  The last inequality can be ensured by fixing a large $n$. Take a large $k$, so
  that the interval
  \[
    \left( |\bar \upsilon_{k n}| {\ln b \over \ln |\lambda_{Z_*,n}| },  |\bar \upsilon_{k n}| {\ln \tilde  b \over \ln |\lambda_{Z_*,n}| }\right)
  \]
  contains an even integer. Now, $l$ can be chosen so that the right endpoint of
  the interval lies between $2 l$ and $2 l+3$.
  
  Next, consider the action of the $k$-th renormalization on the tip $\tau$:
  \begin{align*}
    c_k &= \cRG^k Z(\tau)=A_{k n} (\tau), \\
    \tilde c_k &=\cRG^k \tilde Z(\tau)=\tilde A_{k n}(\tau),
  \end{align*}
  and denote
  \begin{align*}
    X_{k} &= |\pi_1 \tau -\pi_1  c_{k}|, \\
    Y_{k} &= |\pi_2 \tau -\pi_2  c_{k}|=|\pi_2  c_{k}|.
  \end{align*}
  Clearly, since distances $X_k$ and $Y_k$ are those between images of points
  under one iterate of a map $A_{k n}$, which is close to the embedding
  $\iota(\zeta_*)$ of a one-dimensional pair whose both components are
  non-trivial functions of $x$, we have that $X_k \asymp Y_k$.

  We will consider points at three levels of renormalization $\cRG$: level
  $l+k$, level $k$ and level $0$.

  Let $\Phi^{l}_{\cRG^k Z }$ be the map defined in \eqref{Phi}. Consider the
  following points in $\Sigma_{\cRG^k Z}$ and $\Sigma_{\cRG^k \tilde Z}$ at
  level $k$,
  \[
    \tau=\Phi^{l}_{\cRG^k Z}(\tau), \  c_k^{l+k}=\Phi^{l}_{\cRG^k Z}(c_{l+k}) \quad {\rm and} \quad \tilde c_k^{l+k}= \Phi^{l}_{\cRG^k \tilde Z} ( \tilde c_{l+k}).
  \]
  We have by \eqref{lprime} and \eqref{alpha2}, the following inequalities for $k$ large:
  \[
    \left\| \prod_{j=0}^{l-1}   l_{\cRG^{j+k} Z, n}' \circ \prod_{i=j+1}^{l-1} \circ \  l_{\cRG^{i+k} Z,n} \right\|  \asymp |\lambda_{Z_*,n}|^{3 l} \gg  |\lambda_{Z_*,n}|^{3 (l+k)} \asymp C \alpha_2^{l} \gg   \| E_{l,\cRG^k Z} \|.
  \]
Therefore, using \eqref{DPsi}, \eqref{relations}, \eqref{newE} and \eqref{newcE},
  \begin{align}
    \label{dist1} |\pi_1 \tau -\pi_1 c_k^{l+k}| &> \left| c  b^{(1 + O(\alpha^{n}))\absv{\bar{v}_{kn}}} \gamma_1^l Y_{k+l}  -R \prod_{i=0}^{n-1} \sigma_{\cR^i Z_*}^{l} X_{k+l} - C \eps^{2^{k n}} \alpha_2^{l} \right|  \\
    \nonumber &> K b^{(1 + O(\alpha^{n}))\absv{\bar{v}_{kn}}} \gamma_1^l, \\
    \nonumber |\pi_2 \tau -\pi_2 c_k^{l+k}| &> C  |\lambda_{Z_*,n}|^{l}  Y_{k+l}> K  |\lambda_{Z_*,n}|^{l}.
  \end{align}
  Next, consider $c_k$ and
  $$\zeta_k^{l+k} =A_{k n}(c_k^{l+k} ).$$
  We notice that
  \begin{equation}
    A_{k n}= G_{kn} \circ B_{kn},
  \end{equation}
  and since $A_{k n}$ converges to the first map in the pair $Z_*$, given by
  \begin{align*}
    Z_*(x,y) =
    \Lambda\left(  \begin{pmatrix}
        (\eta_*^{r_0}\circ\xi_*)^{r_1}\circ\eta_*(x) \\
        (\eta_*^{r_0}\circ\xi_*)^{r_1-1}\circ\eta_*(x)
      \end{pmatrix}, \quad
    \begin{pmatrix}
      \eta_*^{r_0}\circ\xi_*(x) \\
      x
    \end{pmatrix} \right).
  \end{align*}
  Therefore,
  \begin{equation}
    G_{k n}(x,y)= \left( {g_{k n}(x) + \delta_1(x,y) \atop g_{k n}(y) + \delta_2(x,y)}    \right),
  \end{equation}
  where $g_{k n}$ converges geometrically fast to
  $(\eta_*^{r_0}\circ\xi_*)^{r_1-1}\circ\eta_*(x)$, and the norms of $\delta_1$
  and $\delta_2$ are of order $\epsilon^{2^{k n}}$. We, therefore, obtain that
  \begin{align*}
    \partial_x \pi_2 A_{k n}&= \nabla \pi_2 G_{k n} \circ B_{k n} \cdot \left(\partial_x \pi_1 B_{k n}, \partial_x \pi_2 B_{k n}  \right) \\
    &= \left(\epsilon^{2^{k n}}, O(1)\right) \cdot \left( O(1),O(1)  \right) \\
    &= O(1), \\
    \partial_y \pi_2 A_{k n}&= \nabla \pi_2 G_{k n} \circ B_{k n} \cdot \left(\partial_y \pi_1 B_{k n}, \partial_y \pi_2 B_{k n}  \right) \\
    &= \left(\epsilon^{2^{k n}}, O(1)  \right) \cdot \left( b^{(1+O(\alpha^n))|\bar \upsilon_{k n}| } ,0  \right) \\
    &= \epsilon^{2^{k n}}  b^{(1+O(\alpha^n))|\bar \upsilon_{k n}|}.
  \end{align*}
  Therefore,
  \begin{equation*}
    |\pi_2 c_k -\pi_2 \zeta_k^{l+k}| > K b^{(1 + O(\alpha^{n}))\absv{\bar{v}_{kn}}} \gamma_1^l- C  |\lambda_{Z_*,n}|^{l} b^{(1+O(\alpha^n)) \absv{\bar{v}_{kn}}} \epsilon^{2^{k n}}.
  \end{equation*}
  Since $\gamma_1 \asymp \lambda_{Z_*,n}$ by \eqref{gamma1}, the first term
  dominates the second for large $l$, i.e.
  \begin{equation*}
    |\pi_2 c_k -\pi_2 \zeta_k^{l+k}| > D b^{(1 + O(\alpha^{n}))\absv{\bar{v}_{kn}}} \gamma_1^l.
  \end{equation*}
  Finally, consider
  \[
    z_k^{l+k} =\Phi^{k}_{ Z} (\zeta_k^{l+k}), \quad \xi_k =\Phi^{k}_{Z} (c_k).
  \]
  By \ref{DPsi},
  \begin{equation}\label{conclusion1}
    |z_k^{l+k} - \xi_k  | > c  b^{(1 + O(\alpha^{n}))\absv{\bar{v}_{kn}}} \gamma_1^l \prod_{i=0}^{k-1}  |\lambda_{Z_{i n},n}|    > K b^{(1 + O(\alpha^{n}))\absv{\bar{v}_{kn}}}   |\lambda_{Z_*,n}|^{k} \gamma_1^l
  \end{equation}
  We will now consider similar three pairs of points for the map $\tilde Z$.
  Crucially, the order of the distance between the points $\tau$ and $\tilde
  c_k^{l+k}$ is is different from that for $Z$. By \eqref{relations},
  \begin{align}
   \nonumber |\pi_1 \tau -\pi_1 \tilde c_k^{l+k}| &\le  \left| C  \tilde b^{(1 + O(\alpha^{n}))\absv{\bar{v}_{kn}}} \gamma_1^l \tilde Y_{k+l} +C \eps^{2^{k n}} \alpha_2^{l} \tilde X_{k+l} +R \prod_{k=0}^{n-1} |\sigma_{\cR^i Z_*}|^{l} \tilde X_{k+l} \right|  \\
  \nonumber   & \le  K \prod_{k=0}^{n-1} |\sigma_{\cR^i Z_*}|^{l}   \\
  \label{tdist1}  &\le D  |\lambda_{Z_*,n}|^{3l},\\
  \nonumber |\pi_2 \tau -\pi_2 \tilde c_k^{l+k}| &\le  K  |\lambda_{Z_*,n}|^{l}.
  \end{align}
  Next,
  \begin{align*}
    \partial_x \pi_2 \tilde A_{k n}&= O(1), \\
    \partial_y \pi_2 \tilde A_{k n}&= \epsilon^{2^{k n}}  \tilde b^{(1+O(\alpha^n))|\bar \upsilon_{k n}|},
  \end{align*}
  and
  \begin{align*}
    \partial_x \pi_1 \tilde A_{k n}&= \nabla \pi_1 \tilde G_{k n} \circ \tilde B_{k n} \cdot \left(\partial_x \pi_1 \tilde B_{k n}, \partial_x \pi_2 \tilde B_{k n}  \right) \\
    &= \left(O(1), \epsilon^{2^{k n}} \right) \cdot \left( O(1),O(1)  \right) \\
    &= O(1), \\
    \partial_y \pi_1 \tilde A_{k n}&= \nabla \pi_1 \tilde G_{k n} \circ \tilde B_{k n} \cdot \left(\partial_y \pi_1 \tilde B_{k n}, \partial_y \pi_2 \tilde B_{k n}  \right) \\
    &= \left( O(1) ,\epsilon^{2^{k n}}\right) \cdot \left( \tilde b^{(1+O(\alpha^n))|\bar \upsilon_{k n}| } ,0  \right) \\
    &= \tilde  b^{(1+O(\alpha^n))|\bar \upsilon_{k n}|}.
  \end{align*}
  Therefore, by \eqref{relations},
  \begin{align*}
    |\pi_1 \tilde c_k -\pi_1\tilde \zeta_k^{l+k}| &\le D  |\lambda_{Z_*,n}|^{3l}+   K  \tilde b^{(1+O(\alpha^n))|\bar \upsilon_{k n}|}  |\lambda_{Z_*,n}|^{l} \le K  |\lambda_{Z_*,n}|^{3l},\\
    |\pi_2 \tilde c_k -\pi_2 \tilde \zeta_k^{l+k}| &\le D  |\lambda_{Z_*,n}|^{3l}+  K  \epsilon^{2^{k n}}  \tilde b^{(1+O(\alpha^n))|\bar \upsilon_{k n}|}  |\lambda_{Z_*,n}|^{l} \le K  |\lambda_{Z_*,n}|^{3l}.
  \end{align*}
  Finally,
  \begin{align*}
    |\pi_1 \tilde z_k^{l+k} -\pi_1 \tilde \xi_k  |& < K  |\lambda_{Z_*,n}|^{3l+3 k} +K \eps^{2^{k n}} |\lambda_{Z_*,n}|^{3l}+ C  \tilde b^{(1+O(\alpha^n))|\bar \upsilon_{k n}|} |\lambda_{Z_*,n}|^{3l} \gamma_2^k,  \\
    |\pi_2 \tilde z_k^{l+k} -\pi_2 \tilde \xi_k  |& < K  |\lambda_{Z_*,n}|^{3l+k},
  \end{align*}
  and, by \eqref{relations},
  \begin{equation}\label{conclusion2}
    | \tilde z_k^{l+k} - \tilde \xi_k  | \le  K  |\lambda_{Z_*,n}|^{3l+k},
  \end{equation}
  which, again, by the last inequality of \eqref{relations} implies that
  distances \eqref{conclusion1} and \eqref{conclusion2} between two identical
  pairs points in the attractors $\Sigma$ and $\tilde \Sigma$ are not
  commensurate.

  To complete the proof, assume that the conjugacy $\phi$ is $\kappa$-H\"older.
  Then, necessarily,
  \[
    \absv{z_{k}^{l+k} - \xi_{k}} \le C \absv{\tilde{z}_{k}^{l+k} -
      \tilde{\xi}_{k}}^{\kappa}.
  \]
  for some $C$. By \eqref{conclusion1} and \eqref{conclusion2} this implies
  that
  \begin{equation}
    \label{ineq3}  b^{\absv{\bar{v}_{kn}}}|\lambda_{Z_*,n}|^{k}\gamma_1^{l} \le C |\lambda_{Z_*,n}|^{\kappa (3l+k)}.
  \end{equation}
  Since $\gamma_1 \asymp |\lambda_{Z_*,n}|$, \eqref{ineq3} is equivalent to
  \begin{equation}
    \label{ineq4}  b^{\absv{\bar{v}_{kn}}} |\lambda_{Z_*,n}|^{l+k} \le C |\lambda_{Z_*,n}|^{\kappa (3l+k)}.
  \end{equation}
  We will now bound powers of $|\lambda_{Z_*,n}|$ in \eqref{ineq4} from below
  and above in an optimal way. By \eqref{relations2} we can bound the powers of
  $|\lambda_{Z_*,n}|$ on the right hand side of \eqref{ineq4} using the first
  inequality in \eqref{relations2}:
  \[
    |\lambda_{Z_*,n}|^{3 l} \le C \tilde b^{{3 \over 2}\absv{\bar{v}_{kn}}}.
  \]
  The power $|\lambda_{Z_*,n}|^{l}$ on the left hand side will be bounded from
  below through the second inequality in \eqref{relations2}:
  \[
    |\lambda_{Z_*,n}|^{l} \ge \tilde b^{{1 \over 2} \absv{\bar{v}_{kn}}}.
  \]
  Therefore, \eqref{ineq4} can be weakened as
  \begin{equation}
    \label{ineq5}  b^{\absv{\bar{v}_{kn}}}  \tilde b^{{1 \over 2} \absv{\bar{v}_{kn}}} \le C |\lambda_{Z_*,n}|^{(\kappa-1) k   }    \tilde b^{{3 \over 2} \kappa \absv{\bar{v}_{kn}}},
  \end{equation}
  or, equivalently,
  \begin{equation}
    \label{ineq6}  \absv{\bar{v}_{kn}} \left( \ln b + {1 \over 2} \ln \tilde b -{3 \over 2}  \kappa \ln \tilde b   \right) +(1-\kappa) k \ln |\lambda_{Z_*,n}| \le C.
  \end{equation}
  For a fixed $n$, the first term $\absv{\bar{v}_{kn}} \left( \ln b + ({1 / 2})
    \ln \tilde b -({3 / 2}) \kappa \ln \tilde b \right)$ in $(\ref{ineq6})$
  grows faster in the absolute value than the second term as $k \rightarrow
  \infty$. For this term not to grow arbitrarily large, it is necessary that
  \[
    \ln b + {1 \over 2} \ln \tilde b- {3 \over 2} \kappa \ln \tilde b  \le 0,
  \]
  which implies the claim.
\end{proof}

\appendix

\section{Proof of Proposition \ref*{prop:2Dprojection}}
\label{sec:2DprojectionProof}

Using $q(x) = B(x,0)$ and $p(x) = A(x,0)$ we have
\begin{align*}
  \mathcal{L}([\tilde{A},\tilde{B}])(x) &  = a_{x}(q(x)) \\
                                        & \quad - b(p(x)) - ca_{0}(x) - da_{0}(x)^{2} - ea_{0}(x)^{3} - fa_{0}(x)^{4}. \\
  \mathcal{L}([\tilde{A},\tilde{B}])^{(1)}(x) & = \partial_{1}a_{x}(q(x))\partial_{1}q(x) + \partial_{2}a_{x}(q(x)) \\
                                        & \quad - \partial_{1}b(p(x))a_{0}^{(1)}(x) - \partial_{2}b(p(x))h_{0}^{(1)}(x) \\
                                        & \quad - a_{0}^{(1)}(x)
                                          \left(
                                          c + 2da_{0}(x) + 3ea_{0}(x)^{2} + 4fa_{0}(x)^{3}
                                          \right) \\
  \mathcal{L}([\tilde{A}, \tilde{B}])^{(2)}(x) & = \partial_{1}^{2}a_{x}(q(x))(\partial_{1}q(x))^{2} + 2\partial_{2}\partial_{1}a_{x}(q(x))\partial_{1}q(x) + \partial_{1}a_{x}(q(x))\partial_{1}^{2}q(x) \\
                                        & \quad + \partial_{2}^{2}a_{x}(q(x)) - \left[
                                          \partial_{1}^{2}b(p(x))a_{0}^{(1)}(x)^{2} + 2\partial_{2}\partial_{1}b(p(x))a_{0}^{(1)}(x)h_{0}^{(1)}(x) \right. \\
                                        & \quad + \partial_{2}^{2}b(p(x))h_{0}^{(1)}(x)^{2} + \partial_{1}b(p(x))a_{0}^{(2)}(x) + \partial_{2}b(p(x))h_{0}^{(2)}(x) \\
                                        & \quad + a_{0}^{(2)}(x)
                                          \left(
                                          c + 2da_{0}(x) + 3ea_{0}(x)^{2} + 4fa_{0}(x)^{3}
                                          \right) \\
                                        & \quad \left. + a_{0}^{(1)}(x)^{2}
                                          \left(
                                          2d + 6ea_{0}(x) + 12fa_{0}(x)^{2}
                                          \right)
                                          \right] \\
  \mathcal{L}([\tilde{A}, \tilde{B}])^{(3)}(x) & = \partial_{1}^{3}a_{x}(q(x))q^{(1)}(x)^{3} + 3\partial_{2}\partial_{1}^{2}a_{x}(q(x))q^{(1)}(x)^{2} + 3\partial_{2}^{2}\partial_{1}a_{x}(q(x))q^{(1)}(x) \\
                                        & \quad + \partial_{2}^{3}a_{x}(q(x)) + 3\partial_{1}^{2}a_{x}(q(x))q^{(1)}(x)q^{(2)}(x) \\
                                        & \quad + 3\partial_{2}\partial_{1}a_{x}(q(q))q^{(2)}(x) + \partial_{1}a_{x}(q(x))q^{(3)}(x) \\
                                        & \quad - \left[ \partial_{1}^{3}b(p(x))a_{0}^{(1)}(x)^{3} \right. + 3\partial_{2}\partial_{1}^{2}b(p(x))a_{0}^{(1)}(x)^{2}h_{0}^{(1)}(x) \\
                                        & \quad + 3\partial_{2}^{2}\partial_{1}b(p(x))a_{0}^{(1)}(x)h_{0}^{(1)}(x)^{2} + \partial_{2}^{3}b(p(x))h_{0}^{(1)}(x)^{3} \\
                                        & \quad + 3\partial_{1}^{2}b(p(x))a_{0}^{(1)}(x)a_{0}^{(2)}(x) + 3\partial_{2}\partial_{1}b(p(x))a_{0}^{(2)}(x)h_{0}^{(1)}(x) \\
                                        & \quad + 3\partial_{2}\partial_{1}b(p(x))a_{0}^{(1)}(x)h_{0}^{(2)}(x) + 3\partial_{2}^{2}b(p(x))h_{0}^{(1)}(x)h_{0}^{(2)}(x) \\
                                        & \quad + \partial_{1}b(p(x))a_{0}^{(3)}(x) + \partial_{2}b(p(x))h_{0}^{(3)}(x) \\
                                        & \quad + a_{0}^{(3)}(x)
                                          \left(
                                          c + 2da_{0}(x) + 3ea_{0}(x)^{2} + 4fa_{0}(x)^{3}
                                          \right) \\
                                        & \quad + 3a_{0}^{(1)}(x)a_{0}^{(2)}(x)
                                          \left(
                                          2d + 6ea_{0}(x) + 12fa_{0}(x)^{2}
                                          \right) \\
                                        & \quad + \left. a_{0}^{(1)}(x)^{3}
                                          \left(
                                          6e + 24fa_{0}(x)
                                          \right) \right]
\end{align*}
which gives us:
\begin{align*}
  \mathcal{L}([\tilde{A},\tilde{B}])(0) &  = a_{0}(b_{0}(0)) - b(p(0)) - ca_{0}(0) - da_{0}(0)^{2} - ea_{0}(0)^{3} - fa_{0}(0)^{4}. \\
  \mathcal{L}([\tilde{A},\tilde{B}])^{(1)}(0) & = a_{0}^{(1)}(b_{0}(0))
                                                \left(
                                                b_{0}^{(1)}(0) + c
                                                \right)
                                                + \partial_{2}a_{0}(b_{0}(0)) \\
                                        & \quad - \partial_{1}b(p(0))a_{0}^{(1)}(0) - \partial_{2}b(p(0))h_{0}^{(1)}(0) \\
                                        & \quad - a_{0}^{(1)}(0)
                                          \left(
                                          c + 2da_{0}(0) + 3ea_{0}(0)^{2} + 4fa_{0}(0)^{3}
                                          \right) \\
  \mathcal{L}([\tilde{A}, \tilde{B}])^{(2)}(0) & = a_{0}^{(2)}(b_{0}(0))
                                                 \left(
                                                 b_{0}^{(1)}(0) + c
                                                 \right)^{2}
                                                 + 2\partial_{2}a_{0}^{(1)}(b_{0}(0))
                                                 \left(
                                                 b_{0}^{(1)}(0) + c
                                                 \right) \\
                                        & \quad + a_{0}^{(1)}(b_{0}(0))
                                          \left(
                                          b_{0}^{(2)}(0) + 2d
                                          \right)
                                          + \partial_{2}^{2}a(b_{0}(0)) \\
                                        & \quad - \left[
                                          \partial_{1}^{2}b(p(0))a_{0}^{(1)}(0)^{2} + 2\partial_{2}\partial_{1}b(p(0))a_{0}^{(1)}(0)h_{0}^{(1)}(0) \right. \\
                                        & \quad + \partial_{2}^{2}b(p(0))h_{0}^{(1)}(0)^{2} + \partial_{1}b(p(0))a_{0}^{(2)}(0) + \partial_{2}b(p(0))h_{0}^{(2)}(0) \\
                                        & \quad + a_{0}^{(2)}(0)
                                          \left(
                                          c + 2da_{0}(0) + 3ea_{0}(0)^{2} + 4fa_{0}(0)^{3}
                                          \right) \\
                                        & \quad \left. + a_{0}^{(1)}(0)^{2}
                                          \left(
                                          2d + 6ea_{0}(0) + 12fa_{0}(0)^{2}
                                          \right)
                                          \right] \\
  \mathcal{L}([\tilde{A}, \tilde{B}])^{(3)}(0) & = a_{0}^{(3)}(b_{0}(0))(b_{0}^{(1)}(0) + c)^{3} + 3\partial_{2}a_{0}^{(2)}(b_{0}(0))(b_{0}^{(1)}(0) + c)^{2} \\
                                        & \quad + 3\partial_{2}^{2}a_{0}^{(1)}(b_{0}(0))(b_{0}^{(1)}(0) + c) + \partial_{2}^{3}a_{0}(b_{0}(0)) \\
                                        & \quad + 3a_{0}^{(2)}(b_{0}(0))(b_{0}^{(1)}(0) + c)(b_{0}^{(2)} + 2d) + 3\partial_{2}a_{0}^{(1)}(b_{0}(0))(b_{0}^{(2)} + 2d) \\
                                        & \quad + a_{0}^{(1)}(b_{0}(0))(b_{0}^{(3)}(0) + 6e) - \left[ \partial_{1}^{3}b(p(0))a_{0}^{(1)}(0)^{3} \right. \\
                                        & \quad + 3\partial_{2}\partial_{1}^{2}b(p(0))a_{0}^{(1)}(0)^{2}h_{0}^{(1)}(0) + 3\partial_{2}^{2}\partial_{1}b(p(0))a_{0}^{(1)}(0)h_{0}^{(1)}(0)^{2} \\
                                        & \quad + \partial_{2}^{3}b(p(0))h_{0}^{(1)}(0)^{3} + 3\partial_{1}^{2}b(p(0))a_{0}^{(1)}(0)a_{0}^{(2)}(0) \\
                                        & \quad + 3\partial_{2}\partial_{1}b(p(0))a_{0}^{(2)}(0)h_{0}^{(1)}(0) + 3\partial_{2}\partial_{1}b(p(0))a_{0}^{(1)}(0)h_{0}^{(2)}(0) \\
                                        & \quad + 3\partial_{2}^{2}b(p(0))h_{0}^{(1)}(0)h_{0}^{(2)}(0) + \partial_{1}b(p(0))a_{0}^{(3)}(0) + \partial_{2}b(p(0))h_{0}^{(3)}(0) \\
                                        & \quad + a_{0}^{(3)}(0)
                                          \left(
                                          c + 2da_{0}(0) + 3ea_{0}(0)^{2} + 4fa_{0}(0)^{3}
                                          \right) \\
                                        & \quad + 3a_{0}^{(1)}(0)a_{0}^{(2)}(0)
                                          \left(
                                          2d + 6ea_{0}(0) + 12fa_{0}(0)^{2}
                                          \right) \\
                                        & \quad + \left. a_{0}^{(1)}(0)^{3}
                                          \left(
                                          6e + 24fa_{0}(0)
                                          \right) \right]
\end{align*}

For reference we have:
\begin{align*}
  \mathcal{L}([A,B])(x) & = a_{x}(b_{0}(x)) - b(p(x)) \\
  \mathcal{L}([A,B])^{(1)}(x) & = \partial_{1}a_{x}(b_{0}(x))b_{0}^{(1)}(x) + \partial_{2}a_{x}(b_{0}(x)) \\
                        & - \partial_{1}b(p(x))a_{0}^{(1)}(x) - \partial_{2}b(p(x))h_{0}^{(1)}(x) \\
  \mathcal{L}([A,B])^{(2)}(x) & = \partial_{1}^{2}a_{x}(b_{0}(x))b_{0}^{(1)}(x)^{2} + 2\partial_{2}\partial_{1}a_{x}(b_{0}(x))b_{0}^{(1)}(x) \\
                        & \quad + \partial_{2}^{2}a_{x}(b_{0}(x)) + \partial_{1}a_{x}(b_{0}(x))b_{0}^{(2)}(x) \\
                        & \quad - \left[ \partial_{1}^{2}b(p(x))a_{0}^{(1)}(x)^{2}\right. + 2\partial_{2}\partial_{1}b(p(x))a_{0}^{(1)}(x)h_{0}^{(1)}(x) \\
                        & \quad + \partial_{2}^{2}b(p(x))h_{0}^{(1)}(x)^{2} + \partial_{1}b(p(x))a_{0}^{(2)}(x) \\
                        & \quad + \left. \partial_{2}b(p(x))h_{0}^{(2)}(x) \right] \\
  \mathcal{L}([A,B])^{(3)}(x) & = \partial_{1}^{3}a_{x}(b_{0}(x))b_{0}^{(1)}(x)^{3} + 3\partial_{2}\partial_{1}^{2}a_{x}(b_{0}(x))b_{0}^{(1)}(x)^{2} \\
                        & \quad + 3\partial_{2}^{2}\partial_{1}a_{x}(b_{0}(x))b_{0}^{(1)}(x) + \partial_{2}^{3}a_{x}(b_{0}(x)) + 3\partial_{1}^{2}a_{x}(b_{0}(x))b_{0}^{(1)}(x)b_{0}^{(2)}(x) \\
                        & \quad + 3\partial_{2}\partial_{1}a_{x}(b_{0}(x))b_{0}^{(2)}(x) + \partial_{1}a_{x}(b_{0}(x))b_{0}^{(3)}(x) \\
                        & \quad - \left[ \partial_{1}^{3}b(p(x))a_{0}^{(1)}(x)^{3} \right. + 3\partial_{2}\partial_{1}^{2}b(p(x))a_{0}^{(1)}(x)^{2}h_{0}^{(1)}(x) \\
                        & \quad + 3\partial_{2}^{2}\partial_{1}b(p(x))a_{0}^{(1)}(x)h_{0}^{(1)}(x)^{2} + \partial_{2}^{3}b(p(x))h_{0}^{(1)}(x)^{3} \\
                        & \quad + 3\partial_{1}^{2}b(p(x))a_{0}^{(1)}(x)a_{0}^{(2)}(x) + 3\partial_{2}\partial_{1}b(p(x))a_{0}^{(2)}(x)h_{0}^{(1)}(x) \\
                        & \quad + 3\partial_{2}\partial_{1}b(p(x))a_{0}^{(1)}(x)h_{0}^{(2)}(x) + 3\partial_{2}^{2}b(p(x))h_{0}^{(1)}(x)h_{0}^{(2)}(x) \\
                        & \quad + \partial_{1}b(p(x))a_{0}^{(3)}(x) + \left. \partial_{2}b(p(x))h_{0}^{(3)}(x) \right]
\end{align*}
which gives us:
\begin{align*}
  \mathcal{L}([A,B])(0) & = a_{0}(b_{0}(0)) - b(a_{0}(0), h_{0}(0)) \\
  \mathcal{L}([A,B])^{(1)}(0) & = a_{0}^{(1)}(b_{0}(0))b_{0}^{(1)}(0) + \partial_{2}a_{0}(b_{0}(0)) \\
                        & - \partial_{1}b(a_{0}(0), h_{0}(0))a_{0}^{(1)}(0) - \partial_{2}b(a_{0}(0), h_{0}(0))h_{0}^{(1)}(0) \\
  \mathcal{L}([A,B])^{(2)}(0) & = a_{0}^{(2)}(b_{0}(0))b_{0}^{(1)}(0)^{2} + 2\partial_{2}a_{0}^{(1)}(b_{0}(0))b_{0}^{(1)}(0) \\
                        & \quad + \partial_{2}^{2}a_{0}(b_{0}(0)) + a_{0}^{(1)}(b_{0}(0))b_{0}^{(2)}(0) \\
  & \quad - \left[ \partial_{1}^{2}b(a_{0}(0), h_{0}(0))a_{0}^{(1)}(0)^{2} + 2\partial_{2}\partial_{1}b(a_{0}(0), h_{0}(0))a_{0}^{(1)}(0)h_{0}^{(1)}(0) \right. \\
  & \quad + \partial_{2}^{2}b(a_{0}(0), h_{0}(0))h_{0}^{(1)}(0)^{2} + \partial_{1}b(a_{0}(0), h_{0}(0))a_{0}^{(2)}(0) \\
  & \quad + \left. \partial_{2}b(a_{0}(0), h_{0}(0))h_{0}^{(2)}(0) \right] \\
  \mathcal{L}([A,B])^{(3)}(0) & = a_{0}^{(3)}(b_{0}(0))b_{0}^{(1)}(0)^{3} + 3\partial_{2}a_{0}^{(2)}(b_{0}(0))b_{0}^{(1)}(0)^{2} \\
                        & \quad + 3\partial_{2}^{2}a_{0}^{(1)}(b_{0}(0))b_{0}^{(1)}(0) + \partial_{2}^{3}a_{0}(b_{0}(0)) + 3a_{0}^{(2)}(b_{0}(0))b_{0}^{(1)}(0)b_{0}^{(2)}(0) \\
                        & \quad + 3\partial_{2}a_{0}^{(1)}(b_{0}(0))b_{0}^{(2)}(0) + a_{0}^{(1)}(b_{0}(0))b_{0}^{(3)}(0) \\
                        & \quad - \left[ \partial_{1}^{3}b(p(0))a_{0}^{(1)}(0)^{3} \right. + 3\partial_{2}\partial_{1}^{2}b(p(0))a_{0}^{(1)}(0)^{2}h_{0}^{(1)}(0) \\
                        & \quad + 3\partial_{2}^{2}\partial_{1}b(p(0))a_{0}^{(1)}(0)h_{0}^{(1)}(0)^{2} + \partial_{2}^{3}b(p(0))h_{0}^{(1)}(0)^{3} \\
                        & \quad + 3\partial_{1}^{2}b(p(0))a_{0}^{(1)}(0)a_{0}^{(2)}(0) + 3\partial_{2}\partial_{1}b(p(0))a_{0}^{(2)}(0)h_{0}^{(1)}(0) \\
                        & \quad + 3\partial_{2}\partial_{1}b(p(0))a_{0}^{(1)}(0)h_{0}^{(2)}(0) + 3\partial_{2}^{2}b(p(0))h_{0}^{(1)}(0)h_{0}^{(2)}(0) \\
                        & \quad + \partial_{1}b(p(0))a_{0}^{(3)}(0) + \left. \partial_{2}b(p(0))h_{0}^{(3)}(0) \right]
\end{align*}

The commutation condition of \hyperref[commutation]{Equation \ref*{commutation}}
can then be translated to the following conditions:
\begin{align*}
  0 = \mathcal{L}([\tilde{A},\tilde{B}])(0) & = \mathcal{L}([A,B])(0) - ca_{0}(0) - da_{0}(0)^{2} - ea_{0}(0)^{3} - fa_{0}(0)^{4} \\
  0 = \mathcal{L}([\tilde{A},\tilde{B}])^{(1)}(0) & = \mathcal{L}([A,B])^{(1)}(0) + ca^{(1)}_{0}(b_{0}(0)) \\
                                            & \quad - a_{0}^{(1)}(0)
                                              \left(
                                              c + 2da_{0}(0) + 3ea_{0}(0)^{2} + 4fa_{0}(0)^{3}
                                              \right) \\
  0 = \mathcal{L}([\tilde{A},\tilde{B}])^{(2)}(0) & = \mathcal{L}([A,B])^{(2)}(0) + 2ca_{0}^{(2)}(b_{0}(0))b_{0}^{(1)}(0) + c^{2}a_{0}^{(2)}(b_{0}(0)) \\
                                            & \quad + 2c\partial_{2}a_{0}^{(1)}(b_{0}(0)) + 2da_{0}^{(1)}(b_{0}(0)) \\
                                            & \quad - a_{0}^{(2)}(0)
                                              \left(
                                              c + 2da_{0}(0) + 3ea_{0}(0)^{2} + 4fa_{0}(0)^{3}
                                              \right) \\
                                            & \quad - a_{0}^{(1)}(0)^{2}
                                              \left(
                                              2d + 6ea_{0}(0) + 12fa_{0}(0)^{2}
                                              \right) \\
  0 = \mathcal{L}([\tilde{A},\tilde{B}])^{(3)}(0) & = \mathcal{L}([A,B])^{(3)}(0) + a_{0}^{(3)}(b_{0}(0))(3cb_{0}^{(1)}(0)^{2} + 3c^{2}b_{0}^{(1)}(0) + c^{3}) \\
                                            & \quad + 3\partial_{2}a_{0}^{(2)}(b_{0}(0))(2cb_{0}^{(1)}(0) + c^{2}) + 3c\partial_{2}^{2}a_{0}^{(1)}(b_{0}(0)) \\
                                            & \quad + 3a_{0}^{(2)}(b_{0}(0))(2db_{0}^{(1)}(0) + cb_{0}^{(2)}(0) + 2cd) \\
                                            & \quad + 6d\partial_{2}a_{0}^{(1)}(b_{0}(0)) + 6ea_{0}^{(1)}(b_{0}(0)) \\
                                            & \quad - \left[ a_{0}^{(3)}(0)
                                              \left(
                                              c + 2da_{0}(0) + 3ea_{0}(0)^{2} + 4fa_{0}(0)^{3}
                                              \right) \right. \\
                                            & \quad + 3a_{0}^{(1)}(0)a_{0}^{(2)}(0)
                                              \left(
                                              2d + 6ea_{0}(0) + 12fa_{0}(0)^{2}
                                              \right) \\
                                            & \quad + \left. a_{0}^{(1)}(0)^{3}
                                              \left(
                                              6e + 24fa_{0}(0)
                                              \right) \right]
\end{align*}

We can clearly see that all the equations are satisfied if $\mathcal{L}([A,B]) =
o(\norm{x}^{3})$ and $c=d=e=f=0$. This shows that the projection acts as the
identity on the set of almost commuting pairs and proves the existence of
solutions.

To show that the projection is analytic and uniqueness of solutions we will use
the Implicit Function Theorem. Therefore, consider the map
\[
  (c,d,e,f)\mapsto \left( \mathcal{L}([\tilde{A}, \tilde{B}])(0),
    \mathcal{L}([\tilde{A}, \tilde{B}])^{(1)}(0), \mathcal{L}([\tilde{A},
    \tilde{B}])^{(2)}(0), \mathcal{L}([\tilde{A}, \tilde{B}])^{(3)}(0) \right)
\]
for a fixed pair
$(\tilde{A}, \tilde{B})$. The differential of this map is then given by
\[
  \begin{pmatrix}
    -a_{0}(0) & -a_{0}(0)^{2} & -a_{0}(0)^{3} & -a_{0}(0)^{4} \\
    a_{0}^{(1)}(b_{0}(0)) - a_{0}^{(1)}(0) & -2a_{0}(0)a_{0}^{(1)}(0) &
    -3a_{0}(0)^{2}a_{0}^{(1)}(0) & -4a_{0}(0)^{3}a_{0}^{1}(0) \\
    A & B & C & D \\
    E & F & G & H
  \end{pmatrix}
\]
where we have
\begin{align*}
  A & = 2\left(c + b_{0}^{(1)}(0)\right)a_{0}^{(2)}(b_{0}(0)) + 2\partial_{2}a_{0}^{(1)}(b_{0}(0)) - a_{0}^{(2)}(0) \\
  B & = 2a_{0}^{(1)}(b_{0}(0)) - 2a_{0}^{(2)}(0)a_{0}(0) - 2a_{0}^{(1)}(0)^{2} \\
  C & = -3a_{0}(0)\left(a_{0}(0)a_{0}^{(2)}(0) + 2a_{0}^{(1)}(0)^{2}\right) \\
  D & = -4a_{0}(0)\left(a_{0}(0)a_{0}^{(2)}(0) + 3 a_{0}^{(1)}(0)^{2}\right) \\
  E & = a_{0}^{(3)}(b_{0}(0))
      \left(
      3b_{0}^{(1)}(0)^{2} + 6cb_{0}^{(1)}(0) + 3c^{2}
      \right) + 3\partial_{2}a_{0}^{(2)}(b_{0}(0))
      \left(
      2b_{0}^{(1)}(0) + 2c
      \right) \\
    & \quad + 3\partial_{2}^{2}a_{0}^{(1)}(b_{0}(0)) + 3a_{0}^{(2)}(b_{0}(0))
      \left(
      b_{0}^{(2)}(0) + 2d
      \right) - a_{0}^{(3)}(0) \\
  F & = 3a_{0}^{(2)}(b_{0}(0))
      \left(
      2b_{0}^{(1)}(0) + 2c
      \right) + 6\partial_{2}a_{0}^{(1)}(b_{0}(0)) - 2a_{0}(0)a_{0}^{(3)}(0) - 6a_{0}^{(1)}(0)a_{0}^{(2)}(0) \\
  G & = 6a_{0}^{(1)}(b_{0}(0)) - 3a_{0}(0)^{2}a_{0}^{(3)}(0) - 18a_{0}(0)a_{0}^{(1)}(0)a_{0}^{(2)}(0) - 6a_{0}^{(1)}(0)^{3} \\
  H & = -4a_{0}(0)^{3}a_{0}^{(3)}(0) - 36a_{0}(0)^{2}a_{0}^{(1)}(0)a_{0}^{(2)}(0) - 24a_{0}(0)a_{0}^{(1)}(0)^{3}
\end{align*}

Letting $c=d=e=f=0$ the first two rows stay the same. The last row simplifies as
follows:
\begin{align*}
  A & = 2a_{0}^{(2)}(b_{0}(0))b_{0}^{(1)}(0) + 2\partial_{2}a_{0}^{(1)}(b_{0}(0)) - a_{0}^{(2)}(0) \\
  B & = 2a_{0}^{(1)}(b_{0}(0)) - 2a_{0}^{(2)}(0)a_{0}(0) - 2a_{0}^{(1)}(0)^{2} \\
  C & = -3a_{0}(0)\left(a_{0}(0)a_{0}^{(2)}(0) + 2a_{0}^{(1)}(0)^{2}\right)
\end{align*}

First note that all entries on the first row are bounded away from zero.

For the second row, the first term $a_{0}^{(1)}(b_{0}(0))$ is bounded away from
0 and the rest of the terms in the second row are $O(\epsilon)$ by the Cauchy
integral formula since they are at least $O(\epsilon)$ close to the embedding
of the 1D renormalized pair (which by definition is an almost commuting pair).
Thus the entry in the first column is bounded away from 0 and the rest are
small.

For the third row we have:
\begin{align*}
  A & = 2a_{0}^{(2)}(b_{0}(0))b_{0}^{(1)}(0) + 2\partial_{2}a_{0}^{(1)}(b_{0}(0)) - a_{0}^{(2)}(0) \\
  B & = 2a_{0}^{(1)}(b_{0}(0)) - 2a_{0}^{(2)}(0)a_{0}(0) - 2a_{0}^{(1)}(0)^{2} \\
  C & = -3a_{0}(0)\left(a_{0}(0)a_{0}^{(2)}(0) + 2a_{0}^{(1)}(0)^{2}\right) \\
  D & = -4a_{0}(0)\left(a_{0}(0)a_{0}^{(2)}(0) + 3 a_{0}^{(1)}(0)^{2}\right) \\
  E & = 3a_{0}^{(3)}(b_{0}(0))b_{0}^{(1)}(0)^{2} + 6\partial_{2}a_{0}^{(2)}(b_{0}(0))b_{0}^{(1)}(0) \\
    & \quad + 3\partial_{2}^{2}a_{0}^{(1)}(b_{0}(0)) + 3a_{0}^{(2)}(b_{0}(0))b_{0}^{(2)}(0) - a_{0}^{(3)}(0) \\
  F & = 6a_{0}^{(2)}(b_{0}(0))b_{0}^{(1)}(0) + 6\partial_{2}a_{0}^{(1)}(b_{0}(0)) - 2a_{0}(0)a_{0}^{(3)}(0) - 6a_{0}^{(1)}(0)a_{0}^{(2)}(0) \\
  G & = 6a_{0}^{(1)}(b_{0}(0)) - 3a_{0}(0)^{2}a_{0}^{(3)}(0) - 18a_{0}(0)a_{0}^{(1)}(0)a_{0}^{(2)}(0) - 6a_{0}^{(1)}(0)^{3} \\
  H & = -4a_{0}(0)^{3}a_{0}^{(3)}(0) - 36a_{0}(0)^{2}a_{0}^{(1)}(0)a_{0}^{(2)}(0) - 24a_{0}(0)a_{0}^{(1)}(0)^{3}
\end{align*}
So the entry in the second column is bounded away from zero and the rest are
small.

For the fourth row we have:
\begin{align*}
  E & = O(\epsilon) \\
  F & = O(\epsilon) \\
  G & = 6a_{0}^{(1)}(b_{0}(0)) + O(\epsilon) \\
  H & = O(\epsilon)
\end{align*}
So the entry in the third column is bounded away from zero and the rest are
small.

It follows that the Jacobian matrix at $(c,d,e,f) = (0,0,0,0)$ is invertible for
any commuting $(A,B)$ for small enough $\epsilon$. Using the Implicit Function
Theorem applied to the map
\[
  \left(
    (A,B),c,d,e,f)\mapsto (\mathcal{L}[\tilde{A}, \tilde{B}](0),
    \mathcal{L}([\tilde{A}, \tilde{B}])^{(1)}(0), \mathcal{L}([\tilde{A},
    \tilde{B}])^{(2)}(0), \mathcal{L}([\tilde{A}, \tilde{B}])^{(3)}(0)
  \right)
\]
finishes the proof.

\bibliographystyle{amsalpha}
\bibliography{references}
\end{document}